%

\magnification=\magstep1   
\input amstex
\UseAMSsymbols
\input pictex
\vsize=23truecm
\NoBlackBoxes
\parindent=20pt
  \font\gross=cmbx10 scaled\magstep1 
   \font\rmk=cmr8      \font\ttk=cmtt8


\def\mod{\operatorname{mod}}

\def\Hom{\operatorname{Hom}}

\def\Ext{\operatorname{Ext}}

\def\bdim{\operatorname{\bold{dim}}}

\def\arr#1#2{\arrow <1.5mm> [0.25,0.75] from #1 to #2}

   
\vglue2truecm

\centerline{\gross Distinguished bases of exceptional modules.}
		   		        \bigskip
\centerline{Claus Michael Ringel}     
		  \bigskip\bigskip
{\narrower\narrower An indecomposable representation $M$ of a
quiver $Q= (Q_0,Q_1)$  
is said to be exceptional provided $\Ext^1(M,M) = 0.$ And it is
called a tree module provided 
one can choose a set $\Cal B$ of bases of the vector spaces
$M_x\ (x\in Q_0)$ such that the 
coefficient quiver $\Gamma(M,\Cal B)$ is a tree quiver; we
call $\Cal B$ a tree basis of $M$. It is known that
exceptional modules are tree modules.  
A tree module usually has many tree bases and
the corresponding coefficient quivers may look quite differently.
The aim of this note is to introduce a class of indecomposable modules
which have a distinguished tree basis,
the ``radiation modules'' (generalizing an inductive construction considered
already by Kinser). 
 For a Dynkin
quiver, nearly all indecomposable representations turn out to be radiation modules, 
the only exception is the maximal indecomposable module in case $\Bbb E_8$. Also, 
the exceptional representations of the generalized
Kronecker quivers are given (via the universal cover) by radiation modules. 
Consequently, with the help of Schofield induction 
one can display all the exceptional modules of an
arbitrary quiver in a nice way. 
\par}  
       \bigskip\bigskip
Let $Q = (Q_0,Q_1)$ be a locally finite quiver. We will consider finite-dimensional
representations of $Q$ (thus $kQ$-modules, where $kQ$ is the path algebra of $Q$).
An indecomposable representation $M$ of $Q$  
is said to be {\it exceptional} provided $\Ext^1(M,M) = 0.$ And it is
called a {\it tree module} provided 
one can choose a set $\Cal B$ of bases $\Cal B_x$ of the vector spaces
$M_x\ (x\in Q_0)$ such that the 
coefficient quiver $\Gamma(M,\Cal B)$  is a tree quiver; in this case, we
call $\Cal B$ a {\it tree basis} of $M$. Let me recall the definition of
the coefficient quiver $\Gamma(M,\Cal B)$ as introduced in [R3], 
it is a quiver whose vertices
and arrows are labeled by elements of $Q_0$ and $Q_1$, respectively.
Its vertex set is the disjoint union of the sets $\Cal B_x$, the elements
of $\Cal B_x$ being labeled by $x$. The arrows of $\Gamma(M,\Cal B)$
are obtained as follows: For an arrow $\alpha\:x\to y$ in $Q_1$ and  
$b\in\Cal B_{x},$ write
$M_\alpha(b) = \sum_{b'\in \Cal B_{y}} c_{b'b}b'$ with
coefficients $c_{b'b} \in k$; there is an arrow $b \to b'$ 
in $\Gamma(M,\Cal B)$ with label $\alpha$
provided $c_{b'b} \neq 0.$

It is known [R3] that
exceptional modules are tree modules. 
But even if we know that a module $M$ is a tree module, it often seems
to be difficult to find directly a tree basis. The aim of this note
is to provide for the exceptional modules an algorithm for obtaining
a tree basis. It turns out that we should start by looking 
at indecomposable representations $M$ with a thin vertex
(a vertex $x$ is said to be {\it thin} for $M$ provided
the vector space $M_x$ is one-dimensional). 
The first modules which we will consider are what we call
the radiation modules, see sections 2 and 3.
They are inductively defined,
in any step one constructs an indecomposable module with a thin vertex.
This generalizes a construction 
introduced by Kinser [K] for rooted tree quivers (using the name
``reduced representations'').  

As we will see in section 4, nearly all indecomposable representations of
a Dynkin quiver are radiation modules, the only exception is the maximal
indecomposable module in case $\Bbb E_8$. The radiation modules which
are our main concern are the preprojective and preinjective representations
of a tree without leaves with bipartite orientation (a leaf of a tree is a
vertex with a single neighbor). These modules are studied in sections 5.
Typical examples of trees without leaves are the $n$-regular trees
with $n\ge 2$, note that they occur as the universal cover
of a generalized Kronecker quiver. Thus, 
using covering theory, one obtains a distinguished tree basis
for any exceptional representation of a generalized Kronecker quiver. With the help
of Schofield induction we then get a nice tree basis for any exceptional module,
see the outline in the last section 6.
	\bigskip
{\bf Acknowledgment.} The author is indebted to Valentin Katter and Ryan Kinser 
for spotting errors in preliminary versions of this paper by drawing
the attention to the examples 3 and 5, respectively. 
    \bigskip\bigskip 
{\bf 1. Some exceptional representations $M$ with a thin vertex.}
     \medskip
     If $M$ is a representation of $Q$ and $\dim M_x = 1$ for
some vertex $x$, we say that $x$ is a {\it thin} vertex for $M$.

Let $Q$ be a quiver with underlying graph $\overline Q$ being a tree. Let $x$ be a vertex
of $Q$ and let $Q^x$ be obtained from $Q$ by deleting the vertex $x$ and all the arrows
involving $x$. 

We say that a family of indecomposable modules
$N(1),\dots, N(t)$ is an {\it exceptional family,} provided the modules are pairwise
non-isomorphic and $\Ext^1(N(i),N(j)) = 0$ for all $i, j.$ The family 
$N(1),\dots, N(t)$ is an {\it orthogonal family} provided $\Hom(N(i),N(j)) = 0$ 
for all $i\neq j.$

	\medskip
{\bf Proposition 1.} (a) {\it Let
$N(1),\dots, N(t)$ be an orthogonal family of modules with endomorphism ring $k$, such that
$N(i)_x = 0$ for all $i$,
and assume that for any
index $i$, there is a neighbor $y(i)$ of $x$ with $n(i) = \dim N(i)_{y(i)} > 0$.
Then there is a module $M$ with endomorphism ring $k$, 
unique up to isomorphism, with thin vertex 
$x$ such that the restriction of $M$ to $Q^x$
is $\bigoplus_{i=1}^t N(i)^{n(i)}.$  If the family $N(1),\dots, N(t)$ is
in addition exceptional, then also $M$ is exceptional.} 

(b) {\it Conversely,
let $M$ be an exceptional representation of $Q$, and let $x$ be a vertex of $Q$ with $M_x \neq 0.$
Then the restriction $M'$ of $M$ to $Q^x$ is of the form $M' = \bigoplus_{i=1}^t N(i)^{n(i)}$, 
where $N(1),\dots, N(t)$ is an exceptional family of indecomposable modules
and $n(i) > 0$ for all $i$.
For any
index $i$, there is a unique neighbor $y(i)$ of $x$ with $\dim N(i)_{y(i)}\neq 0.$
If $\dim M_x = 1$, then $n(i) \le \dim N(i)_{y(i)}$ for all $ 1\le i \le t$. 
If $\dim M_x = 1$ and $\dim N(i)_{y(i)} = 1$ for $i\in I$, 
where $I$ is a subset of $\{1,2,\dots,t\}$,
then the modules $N(i)$ with $i\in I$
form an orthogonal family.}
   \bigskip
Proof: (a) Let 
$N(1),\dots, N(t)$ be an orthogonal family of modules with endomorphismring $k$
such that $N(i)_x = 0$, and assume
that for any
index $i$, there is a neighbor $y(i)$ of $x$ with $n(i) = \dim N(i)_{y(i)} > 0$.
We can assume that there are arrows $y(i) \to x$ for $1 \le i \le s$ and
arrows $x \to y(i)$ for $s+1 \le i \le t$.
Since $N(i)_x = 0$, we see that the modules  $S(x), N(1), \dots, N(t)$ 
are orthogonal bricks, thus we can use the process of simplification 
(see [R1]): it asserts that the subcategory $\Cal A$ of representations
of $Q$ which have a filtration with factors in $\{S(x),N(1),\dots,N(t)\}$
is an exact abelian subcategory of $\mod kQ$ which is closed under
extensions, and which has the simple
objects $S(x), N(1), \dots, N(t)$. 
Since the simple objects of $\Cal A$
have endomorphism ring $k$, the category $\Cal A$ is
equivalent to the category of locally nilpotent representation of a quiver $Q(\Cal A)$
with $t+1$ vertices which are labeled $S(x), N(1),\dots, N(t)$. For every pair
$N,N'$ of vertices $Q(\Cal A)$, the number of arrows from $N$ to $N'$ is equal to
$\dim\Ext^1(N,N')$. For $1\le i \le s$, we have $\dim\Ext^1(N(i),S(x)) = \dim N(i)_{y(i)}
= n(i)$ and $\dim\Ext^1(S(s),N(i)) = 0$. 
For $s+1\le i \le t$, we have $\dim\Ext^1(S(x),N(i)) = \dim N(i)_{y(i)}
= n(i)$ and $\dim\Ext^1(N(i),S(s)) = 0$, thus the quiver $Q(\Cal A)$ is
obtained from the following quiver $Q'$ by adding arrows $N(i) \to N(j)$ for $1\le i,j\le t$,
according to the extension groups $\Ext^1(N(i),N(j))$ in question.
$$
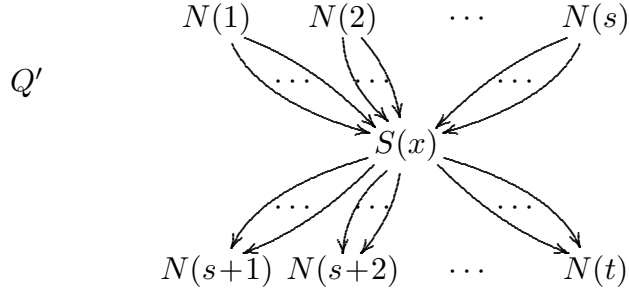
{\beginpicture
\setcoordinatesystem units <.7cm,.7cm>
\put{$N(1)$} at 0 4
\put{$N(2)$} at 2 4
\put{$N(s)$} at 6 4
\put{$N(s\!+\!1)$} at 0 0
\put{$N(s\!+\!2)$} at 2 0
\put{$N(t)$} at 6 0
\put{$S(x)$} at 3 2
\multiput{$\cdots$} at 4 0  4 4 /
\setquadratic
\plot 0.25 3.6  1 2.8  2.4 2.2 /
\plot 0.45 3.7  1.5 3.15  2.5 2.3 /
\plot 2 3.7  2.2 3   2.7 2.4 /
\plot 2.3 3.7  2.7 3.1   2.9 2.45 /

\plot 5.75 3.6  5 2.8  3.6 2.2 /
\plot 5.55 3.7  4.5 3.15  3.5 2.3 /

\arr{2.4 2.2}{2.42 2.195}
\arr{2.5 2.3}{2.51 2.29}
\arr{2.7 2.4}{2.71 2.39}
\arr{2.9 2.45}{2.902 2.44}

\arr{3.5 2.3}{3.48 2.283}
\arr{3.6 2.2}{3.59 2.197}

\plot 0.25 0.4  1 1.2  2.4 1.8 /
\plot 0.45 .3  1.5 .85  2.5 1.7 /
\plot 2 .3  2.2 1   2.7 1.6 /
\plot 2.3 .3  2.7 .9   2.9 1.55 /

\plot 5.75 .4  5 1.2  3.6 1.8 /
\plot 5.55 .3  4.5 .85  3.5 1.7 /

\arr{0.25 0.4}{0.24 0.38}  
\arr{0.45 .3}{0.43 .293}  
\arr{2 .3}{2 .29} 
\arr{2.3 .3}{2.292 .29} 

\arr{5.55 .3}{5.57 .293}  
\arr{5.75 .4}{5.755 .39}

\multiput{$\cdots$} at 1.24 3  2.5 3  4.76 3 
 1.24 1  2.5 1  4.76 1 /
\put{$Q'$} at -3 3 
\endpicture}
$$
Let me repeat that the number of arrows of $Q'$ between $N(i)$ and $S(x)$
is $n(i)$, for $1\le i \le t$.

However, we are only interested in the subcategory $\Cal A'$ of $\Cal A$ consisting
of the representations of $Q$ whose restriction to $Q^x$ is a direct sum of modules
of the form $\bigoplus_i N(i)^{m(i)}$ with $m(i)\in \Bbb N_0$. This subcategory
$\Cal A'$ is an exact abelian subcategory of $\Cal A$ with the same simple 
objects and $\Cal A'$ is 
equivalent to the category of locally nilpotent representation of the 
quiver $Q'$ itself.

Now the quiver $Q'$ has the following real root:
$$
\hbox{\beginpicture
\setcoordinatesystem units <.7cm,.3cm>
\put{$n(1)$} at 0 4
\put{$n(2)$} at 2 4
\put{$n(s)$} at 6 4
\put{$n(s\!+\!1)$} at 0 0
\put{$n(s\!+\!2)$} at 2 0
\put{$n(t)$} at 6 0
\put{$1$} at 3 2
\multiput{$\cdots$} at 4 0  4 4 /
\endpicture} 
$$
thus in $\Cal A'$ there is precisely 
one indecomposable object $M$ with one factor $S(x)$ and $n(i)$ factors
$N(i)$, for $1\le i \le t$ and the endomorphism ring of $M$ is $k$.

In case the family $N(1),\dots,N(t)$ is in addition exceptional, then 
$\Cal A' = \Cal A$ and therefore $Q(\Cal A) = Q'$. In this case, the object 
$M$ is an exceptional module. Namely, the module  $M$ is always 
exceptional when  considered as an object in $\Cal A'$; 
now $\Cal A' = \Cal A$, and $\Cal A$ is closed under extensions in $\mod kQ$.
 
	\medskip
(b) We use the following lemma:
    \medskip
{\bf Lemma 1.} {\it Let $X'$ be a submodule of $X$ and $X'' = X/X'.$ We assume that  $\Ext^1(X,X) = 0$
and $\Hom(X',X'') = 0.$ Then $\Ext^1(X',X') = 0$ and $\Ext^1(X'',X'') = 0.$}
    \medskip
Proof: The inclusion map $X' \to X$ yields an epimorphism
$\Ext^1(X,X) \to \Ext^1(X',X)$, thus $\Ext^1(X',X) = 0.$ Applying $\Hom(X',-)$ to
the exact sequence $0 \to X' \to X \to X'' \to 0$ yields an exact
sequence
$$
 \Hom(X',X'') \to \Ext^1(X',X') \to \Ext^1(X',X).
$$
Since the end terms are zero, also the middle term is zero. This shows that 
$\Ext^1(X',X') = 0.$ By duality, also $\Ext^1(X'',X'') = 0.$
    \medskip
Assume now that $M$ is an exceptional representation of $Q$.
Let $y_1,\dots, y_n$ be the neighbors of $x$, say with arrow $y_i\to x$ for $1\le i \le r$
and $x \to y_i$ for $r+1\le i \le n.$ The quiver $Q^x$ is the disjoint union of two parts
$Q^+$ and $Q^-$, where $Q^-$ is the union of the connected components of $Q^x$ which 
contain the vertices $y_1,\dots,y_r$, and $Q^+$ 
is the union of the connected components of $Q^x$ which 
contain the vertices $y_{r+1},\dots,y_n$. Let $M^+$ be the restriction of $M$ to $Q^+$
and $M^-$ the restriction of $M$ to $Q^-$. Note that $M^+$ is a submodule of $M$, whereas
$M^-$ is a factor module of $M$, say $M^- = M/X$. Thus, $M$ has the
following chain of submodules $M^+ \subseteq X \subseteq M$ and $X/M^+$ is a direct sum
of copies of $S(x)$. Now we apply Lemma 1 to the submodule $X$ of $M$, this is possible,
since $X$ and $M/X$ has disjoint support, thus $\Hom(X,M/X) =0$. We see that
$\Ext^1(X,X) = 0$ and $\Ext^1(M^-,M^-) = 0.$ Next, we apply Lemma 1 to the submodule 
$M^+$ of $X$, again using that we deal with modules $M^+$ and $X/M^+$ with disjoint support.
We conclude that $\Ext^1(M^+,M^+) = 0.$ Since there is no arrow between vertices
of $Q^-$ and $Q^+$, thus $\Ext^1(M^-,M^+) = 0 = \Ext^1(M^+,M^-)$. Note that $M' = M^+\oplus M^-$.
Altogether we have 
shown that $\Ext^1(M',M') = 0,$ thus the modules $N(1),\dots, N(t)$ form
an exceptional family. 

Now consider a module $N(i)$.
Since the support of $M$ is connected, and $M$ is indecomposable, 
at least one of the vertices $y_j$ must belong to the support of $N(i)$. Since $\overline Q$ is a tree, 
there is just one such vertex. This shows that there is a unique neighbor $y(i)$ of
$x$ such that $N(i)_{y(i)} \neq 0$, namely $y(i) = y_j.$ 

From now on, let us assume that $\dim M_x = 1.$ 
Let us show that $n(i) \le \dim N(i)_{y(i)}$ for all $1\le i \le t$. 
We assume that for $1\le i \le s$, the arrow between $x$ and $y(i)$ ends in $x$, whereas
for $s+1\le i \le t$, it starts in $x$. First, consider some $i$ with $s+1\le i \le t$, say
$i = r+1.$ Then $\dim \Ext^1(S(x),N(r+1)) = \dim N(r+1)_{y(r+1)}$. There is the exact sequence
$$
 0 \to N(r+1)^{n(r+1)} \oplus X' \to X \to S(x) \to 0 \quad \text{with}\quad X' =
\bigoplus_{i=r+2}^t N(i)^{n(i)}.
$$
If $n(i) > \dim\Ext^1(S(x),N(r+1)),$ then $X$ splits off a copy of $N(r+1)$ (see the
$\Ext$-Lemma in [R4]). But this is impossible, since $X$ is indecomposable.
The dual argument works for $1\le i \le r.$ This shows that $n(i)\le \dim N(i)_{y(i)}$
for all indices $1\le i\le t$.

In particular, we see that $\dim N(i)_{y(i)} = 1$ implies $n(i) = 1$. 
Now assume in addition that there are two different indices $i,j$ with 
$\dim N(i)_x = 1 = \dim N(j)_x$, say $i =1, j = 2.$ We want to show that $\Hom(N(1),N(2) = 0.$
Thus, assume that $\Hom(N(1),N(2)) \neq 0.$ Then the support of $N(1)$ and $N(2)$ is contained
in the same connected component of $Q^x$, thus $y(1) = y(2)$. Let $y = y(1) = y(2)$.

Let us assume that the arrow $\alpha$ between $x$ and $y$ starts in $x.$ 
Since $\dim M_x = 1,$ we can assume that $M_x = k$ and consider the
element $M_\alpha(1)\in M_y$. Note that 
$$
 M_y = M'_y = \bigoplus_{i=1}^t N(i)_y^{n(i)} = N(1)_y\oplus N(2)_y \oplus N'_y
 \quad\text{with} \quad N' = \bigoplus_{i=2}^t N(i)_y^{n(i)},
$$
thus we can write $M_\alpha(1) = (a_1,a_2,a_3)$ with $a_1 \in N(1)_y,\ a_2\in N(2)_y,\ 
a_3\in N'_y.$ Note that both $a_1,a_2$ are non-zero, since otherwise $M$ would split
off $N(1)$ or $N(2)$, respectively. Let $\phi\:N(1) \to N(2)$
be a non-zero homomorphism. Since $\Ext^1(N(2),N(1)) = 0$, we know that $\phi$ has be a
monomorphism or an epimorphism. It follows from $\dim N(1)_y = 1 = \dim N(2)_y$,
that $\phi_y$ is bijective. Replacing if necessary $\phi$ by a scalar multiple,
we can assume that $\phi_y(a_1) = a_2$. Let $N'' \subseteq N(1)\oplus N(2)$ be the
graph of $\phi$, thus $N''_z = \{(b,\phi(b))\mid b\in N(1)_z\}$ for all vertices $z$
of $Q$. The module $N''$
is a submodule of $N(1)\oplus N(2)$ which has $0\oplus N(2)$ as a direct
complement inside $N(1)\oplus N(2)$. Since we see that $M_\alpha$ maps $M_x$ into
$N''_y\oplus N'_y$, it follows that $M$ splits off $N(2)$, namely $M = Y \oplus N(2)$,
where $Y_x = M_x$ and the restriction of $Y$ to $Q^x$ is equal to $N''\oplus N'.$
This contradicts the indecomposability of $M$.

We use the dual considerations in case 
the arrow between $x$ and $y$ starts in $y$ and
terminates in $x$. This concludes the proof.

   \bigskip
{\bf Example 1.} Consider the following exceptional representation $M$, encircled
is a vertex $x$ with $\dim M_x = 1$, note that $x$ has a unique neighbor $y$.
The lower line exhibits the dimension vectors
of the modules $N(i)$, here the vertex $y$ is encircled
by a dotted circle.
$$
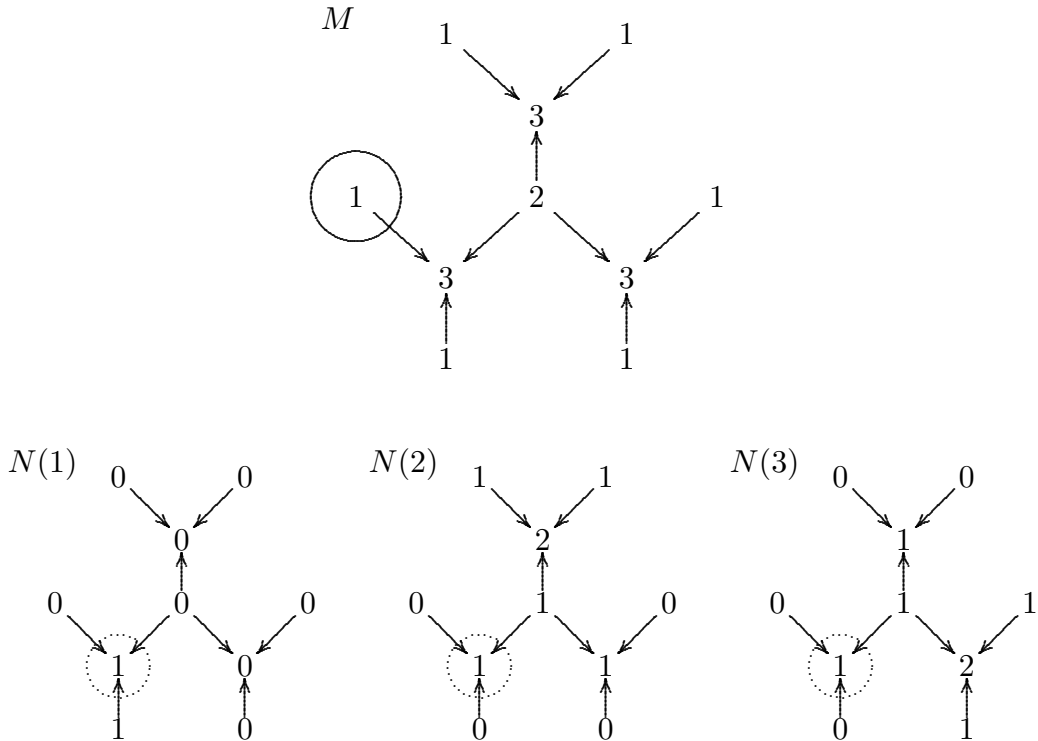
\beginpicture
\setcoordinatesystem units <1cm,.9cm>
\put{\beginpicture
\put{$M$} at -2.2 2.2
\multiput{$1$} at  -1 2  1 2  -2 0  2 0  -1 -2  1 -2 /
\multiput{$2$} at 0 0   /
\multiput{$3$} at 0 1  -1 -1  1 -1   /

\arr{0 0.2}{0 0.8}
\arr{-.2 -.2}{-.8 -.8}
\arr{.2 -.2}{.8 -.8}
\arr{-.8 1.8}{-.2 1.2}
\arr{.8 1.8}{.2 1.2}
\arr{-1.8 -.2}{-1.2 -.8}
\arr{1.8 -.2}{1.2 -.8}
\arr{-1 -1.8}{-1 -1.2}
\arr{1 -1.8}{1 -1.2}
\circulararc 360 degrees from -1.5 0  center at -2 0 
\endpicture} at 0 0 
\put{\beginpicture
\setcoordinatesystem units <.7cm,.7cm>
\put{$N(1)$} at -2.2 2.2
\multiput{$0$} at  0 0  0 1  1 -1   -1 2  1 2  -2 0  2 0  1 -2 /
\multiput{$1$} at    -1 -1   -1 -2  /

\arr{0 0.2}{0 0.8}
\arr{-.2 -.2}{-.8 -.8}
\arr{.2 -.2}{.8 -.8}
\arr{-.8 1.8}{-.2 1.2}
\arr{.8 1.8}{.2 1.2}
\arr{-1.8 -.2}{-1.2 -.8}
\arr{1.8 -.2}{1.2 -.8}
\arr{-1 -1.8}{-1 -1.2}
\arr{1 -1.8}{1 -1.2}
\setdots <.7mm>
\circulararc 360 degrees from -.5 -1  center at -1 -1 
\endpicture} at -4 -5
\put{\beginpicture
\setcoordinatesystem units <.7cm,.7cm>
\put{$N(2)$} at -2.2 2.2
\multiput{$0$} at     -2 0  2 0  -1 -2  1 -2 /
\multiput{$1$} at   -1 -1   0 0  1 -1   -1 2  1 2 /
\put{$2$} at   0 1  
\arr{0 0.2}{0 0.8}
\arr{-.2 -.2}{-.8 -.8}
\arr{.2 -.2}{.8 -.8}
\arr{-.8 1.8}{-.2 1.2}
\arr{.8 1.8}{.2 1.2}
\arr{-1.8 -.2}{-1.2 -.8}
\arr{1.8 -.2}{1.2 -.8}
\arr{-1 -1.8}{-1 -1.2}
\arr{1 -1.8}{1 -1.2}
\setdots <.7mm>
\circulararc 360 degrees from -.5 -1  center at -1 -1 
\endpicture} at 0 -5
\put{\beginpicture
\setcoordinatesystem units <.7cm,.7cm>
\put{$N(3)$} at -2.2 2.2
\multiput{$0$} at     -1 2  1 2  -2 0  -1 -2 /
\multiput{$1$} at  0 0  -1 -1  0 1  1 -2  2 0  /
\put{$2$} at   1 -1

\arr{0 0.2}{0 0.8}
\arr{-.2 -.2}{-.8 -.8}
\arr{.2 -.2}{.8 -.8}
\arr{-.8 1.8}{-.2 1.2}
\arr{.8 1.8}{.2 1.2}
\arr{-1.8 -.2}{-1.2 -.8}
\arr{1.8 -.2}{1.2 -.8}
\arr{-1 -1.8}{-1 -1.2}
\arr{1 -1.8}{1 -1.2}
\setdots <.7mm>
\circulararc 360 degrees from -.5 -1  center at -1 -1 
\endpicture} at 4 -5
\endpicture
$$
(A related representation, namely the representation $P(c,3)$ for the
$3$-regular tree with bipartite orientation and $c$ a source, will be discussed
in detail towards the end of section 5. The representation $M$ considered
above is the restriction of $P(c,3)$ to the ball with center $c$ and radius $2$.)
	\bigskip
We mainly will be interested in exceptional modules $M$ with a thin vertex $x$ such that
the restriction $M'$  to $Q^x$ 
decomposes as the direct sum of indecomposable modules $N(i)$ with $\dim N(i)_{y(i)} = 1,$
so that also $n(i) = 1.$ But let us exhibit here an example with $\dim N(i)_{y(i)} > 1$
for some index $i$ (further examples will be provided in section 3, namely 
part 2 of example 3, as well as example 4 --- whereas the indecomposable direct summands
of $M$ in example 2 are orthogonal, those in section 3 are not).
    \medskip
{\bf Example 2.} An exceptional module $M$ with a thin vertex 
$x$ such that the restriction $M'$ of $M$ to
$Q^x$ is of the form $N(1)\oplus N(2)^2$ for an exceptional
orthogonal family $N(1), N(2)$.
$$
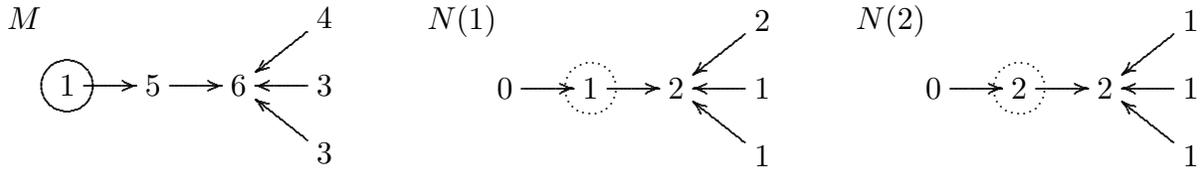
\beginpicture
\setcoordinatesystem units <.95cm,1.5cm>
\put{\beginpicture
\put{$M$} at -.5 0.5
\put{$1$} at 0 0
\put{$5$} at 1 0
\put{$6$} at 2 0
\put{$4$} at 3 0.5
\put{$3$} at 3 0
\put{$3$} at 3 -.5
\arr{0.2 0}{0.8 0}
\arr{1.2 0}{1.8 0}
\arr{2.8 0}{2.2 0}
\arr{2.8 0.4}{2.2 0.1}
\arr{2.8 -.4}{2.2 -.1}
\circulararc 360 degrees from 0.3 0  center at 0 0 
\endpicture} at 0 0
\put{\beginpicture
\put{$N(1)$} at -.5 0.5
\put{$0$} at 0 0
\put{$1$} at 1 0
\put{$2$} at 2 0
\put{$2$} at 3 0.5
\put{$1$} at 3 0
\put{$1$} at 3 -.5
\arr{0.2 0}{0.8 0}
\arr{1.2 0}{1.8 0}
\arr{2.8 0}{2.2 0}
\arr{2.8 0.4}{2.2 0.1}
\arr{2.8 -.4}{2.2 -.1}
\setdots <.7mm>
\circulararc 360 degrees from 1.3 0  center at 1 0 
\endpicture} at 5 0
\put{\beginpicture
\put{$N(2)$} at -.5 0.5
\put{$0$} at 0 0
\put{$2$} at 1 0
\put{$2$} at 2 0
\put{$1$} at 3 0.5
\put{$1$} at 3 0
\put{$1$} at 3 -.5
\arr{0.2 0}{0.8 0}
\arr{1.2 0}{1.8 0}
\arr{2.8 0}{2.2 0}
\arr{2.8 0.4}{2.2 0.1}
\arr{2.8 -.4}{2.2 -.1}
\setdots <.7mm>
\circulararc 360 degrees from 1.3 0  center at 1 0 
\endpicture} at 10 0
\endpicture
$$
	\bigskip\bigskip 
{\bf 2. Radiation modules.}
     \medskip
Let us assume that we deal with a quiver $Q$ which is a tree.
We consider pairs $(M,x)$ where $M$ is an indecomposable module
and $x$ is a vertex with $\dim M_x = 1$ (called the {\it origin} of $M$ or better of the pair
$(M,x)$).
We define the class of 
{\it radiation modules} as well as the corresponding {\it radiation quivers}
inductively as follows. 

First of all, for any vertex $x$, the pair $(S(x),x)$ is a radiation module, its radiation
quiver $R(S(x),x))$ is the quiver with a single vertex (and no arrow), the vertex
being labeled $x$. 
Second, if $M$ is an indecomposable module of length at least 2
and $x$ is a vertex with $\dim M_x = 1,$ 
then $(M,x)$ is a radiation module provided the restriction of $M$ to the
quiver obtained from $Q$ by deleting the vertex $x$ is the direct sum
$\bigoplus N(i)$ of an orthogonal family of indecomposable 
modules $N(i)$ such that for any index
$i$ there is a neighbor $y(i)$ of $x$ with $(N(i),y(i))$ being a
radiation module. The radiation quiver $R(M,x)$ of $(M,x)$ is obtained from the
disjoint union of the radiation quivers $R(N(i),y(i))$ by adding a vertex with
label $x$ and for every index $i$ we connect $x$ and $y(i)\in R(N(i),y(i))_0$
by an arrow which points
in the same direction as the arrow in $Q$ between $x$ and $y(i)$, and
we use this arrow of $Q$ also as the label of the arrow in $R(M,x)$. Note that
the labels of the vertices and the arrows of $R(M,x)$ provide a quiver
homomorphism $R(M,x) \to Q.$ 
	     \medskip
Let $(M,x)$ be a radiation module, and $b$ a non-zero element of $M_x$ (note that $M_x$
is one-dimensional). 
For any radiation module $(M,x)$ we define a basis $\Cal B(M,x)$ 
as follows, we call it a {\it radiation basis} containing $b$.
Assume that 
the restriction of $M$ to the quiver obtained by deleting the vertex $x$ is
$\bigoplus N(i)$ with indecomposable representations $N(i)$. For every index $i$,
there is a vertex $y(i)$ which is a neighbor of $x$ such that $(N(i),y(i))$ is a radiation module. 
Since $x,y(i)$ are neighbors, there is an arrow say $\alpha$
which connects them. If $\alpha\:x \to y(i)$,  
let $\Cal B(i)$ be a radiation basis of $(N(i),y(i))$ containing the element $M_\alpha(b).$ 
If $\alpha\:y(i) \to x$,  let $\Cal B(i)$
 be a radiation basis of $(N(i),y(i))$ containing the element $M_\alpha^{-1}(b).$ 
Finally, let $\Cal B(M,x)$ be the set
containing the element $b$ as well as all the elements in the 
(disjoint!) sets $\Cal B(i)$.
	\medskip 
{\bf Proposition 2.} {\it If $(M,x)$ is a radiation module, then the endomorphism ring of $M$ is 
$k$. 
A radiation basis $\Cal B(M,x)$ of $M$ is a tree basis and the corresponding coefficient quiver
is $R(M,x)$.}
    \medskip
Proof. The first assertion is shown by induction. 
If $M = S(x)$, then the endomorphism ring of $M$ is $k$. If 
$M$ is of length at least 2, then the restriction of $M$ to the
quiver $Q^x$ is the direct sum of radiation modules $(N(i),y(i))$, and by
induction the endomorphism ring of $N(i)$ is $k$. Now the
family of modules $N(i)$ is an orthogonal family of modules with 
endomorphism ring $k$ and for any index
$i$ there is the neighbor $y(i)$ of $x$ with $n(i) = \dim N(i)_{y(i)} = 1.$
Thus $M$ is the module constructed in Proposition 1(a) using the vertex
$x$ and the modules $N(i)$. According to this proposition, 
the endomorphism ring of $M$ is $k$. 

Also the remaining assertions are
straight-forward: Namely, assume that there is given 
a finite set of trees 
$T(i)$ and for every $i$ a vertex $y(i)$ in $T(i)$. If we take the disjoint union of 
the graphs $T(i)$, an additional vertex $x$ and connect $x$ with the vertices $y(i)$
by a single edge, then we obtain again a tree.
In order to see that $R(M,x)$ is a coefficient quiver of $M,$
one observes by induction that $\Cal B(M,x)$ is a basis
of $M$, the corresponding coefficient quiver is just $R(M,x)$.
r is just $R(M,x)$.
	\bigskip
{\bf Remark.} In the special situation of $Q$ being 
a tree quiver with precisely one root, radiation modules 
have been considered before by Kinser [K], he called them
reduced representations. His detailed investigation
of these modules is of great interest. Note that they
play a prominent role in his study
of the representation ring of such a quiver. See also
the paper [KM] by Katter and Mahrt for a further study of 
this class of radiation modules.  
	\bigskip\bigskip
{\bf 3. Exceptional radiation modules.}
	\medskip
We are mainly interested in modules which are both exceptional as well as radiation modules,
thus in exceptional modules with a thin vertex. We should stress that by construction radiation 
modules have a thin vertex, but may not be exceptional, and conversely, there are exceptional
modules with a thin vertex which are not radiation modules. Here are corresponding examples:
	\medskip 
{\bf Example 3.} An exceptional module $M$ with a thin vertex, 
such that $(M,x)$ is not a radiation module for any thin vertex $x$. The module $M$
which we exhibit is an indecomposable preinjective module with a unique thin vertex.
On the right, we
show the indecomposable direct summands $N(1),N(2)$ of $M'$, here, the
neighbor $y = y(i)$ of $x$ is encircled using a dotted circle. 
Note that the modules $N(1),N(2)$ form an exceptional family. Since 
$\Hom(N(1),N(2)) \neq 0,$ the family is not orthogonal.
$$
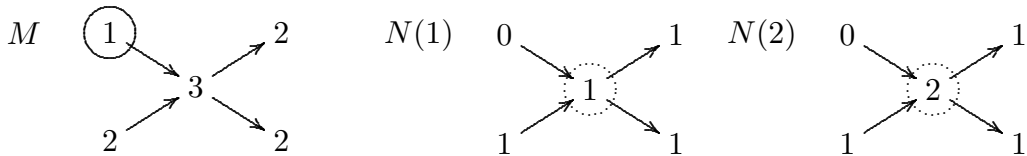
\beginpicture
\setcoordinatesystem units <.95cm,.6cm>
\put{\beginpicture
\put{$M$} at -1 1
\put{$1$} at 0 1
\put{$2$} at 0 -1
\put{$3$} at 1 0
\put{$2$} at 2 1
\put{$2$} at 2 -1
\arr{.2 0.8}{.8 0.2}
\arr{.2 -.8}{.8 -.2}
\arr{1.2 0.2}{1.8 0.8}
\arr{1.2 -.2}{1.8 -.8}

\circulararc 360 degrees from .3 1  center at 0 1 
\endpicture} at 0 2
\put{\beginpicture
\put{$N(1)$} at -1 1
\put{$0$} at 0 1
\put{$1$} at 0 -1
\put{$1$} at 1 0
\put{$1$} at 2 1
\put{$1$} at 2 -1
\arr{.2 0.8}{.8 0.2}
\arr{.2 -.8}{.8 -.2}
\arr{1.2 0.2}{1.8 0.8}
\arr{1.2 -.2}{1.8 -.8}
\setdots <.7mm>
\circulararc 360 degrees from 1.3 0  center at 1 0 
\endpicture} at 4.5 2
\put{\beginpicture
\put{$N(2)$} at -1 1
\put{$0$} at 0 1
\put{$1$} at 0 -1
\put{$2$} at 1 0
\put{$1$} at 2 1
\put{$1$} at 2 -1
\arr{.2 0.8}{.8 0.2}
\arr{.2 -.8}{.8 -.2}
\arr{1.2 0.2}{1.8 0.8}
\arr{1.2 -.2}{1.8 -.8}
\setdots <.7mm>
\circulararc 360 degrees from 1.3 0  center at 1 0 
\endpicture} at 8.5 2

\endpicture
$$
Also, one may look at this example in the light of Proposition 1(b): 
If we write $M' = N(1)^{n(1)}\oplus N(2)^{n(2)}$, then $n(1) = n(2) = 1$. Now 
$\dim N(2)_{y(i)} = 2$, thus $n(2) < \dim N(2)_{y(i)}.$ 
	\medskip 
{\bf Example 4.} This is an example of an exceptional module $M$ 
with vertices $x,x'$ such that
$(M,x)$ is a radiation module, whereas $(M,x')$ is not. The module $M$ which
we present here is indecomposable and preinjective, the vertices $x,x'$ are
encircled. We write $M^x$ for the restriction of $M$ to
$Q^x$, and $M^{x'}$ for its restriction to $Q^{x'}$. On the right, we
show the indecomposable direct summands $N(i)$ and $N'(i)$ of $M^x$ and
$M^{x'}$ respectively. Any of the vertices $x,x'$ has a unique neighbor $y$,
this vertex is encircled using a dotted circle. 
Note that the modules $N(1),N(2),N(3)$ form an
orthogonal exceptional family, whereas $N'(1), N'(2)$ is an exceptional
family which is not orthogonal (namely, we have $\dim\Hom(N'(1),N'(2)) = 1$). 
$$
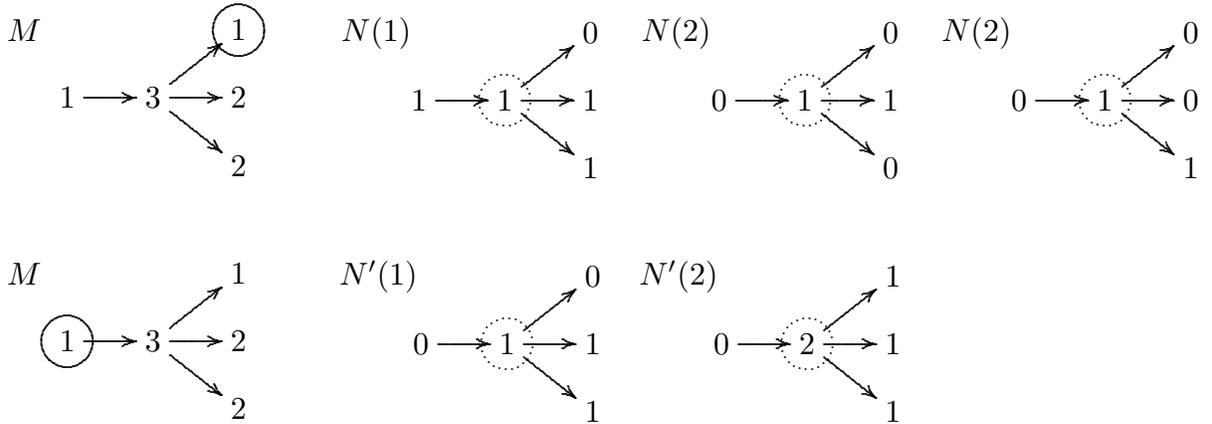
\beginpicture
\setcoordinatesystem units <.95cm,1.5cm>
\put{\beginpicture
\put{$M$} at .5 0.5
\put{$1$} at 1 0
\put{$3$} at 2 0
\put{$1$} at 3 0.5
\put{$2$} at 3 0
\put{$2$} at 3 -.5
\arr{1.2 0}{1.8 0}
\arr{2.2 0}{2.8 0}
\arr{2.2 0.1}{2.8 0.4}
\arr{2.2 -.1}{2.8 -.4}
\circulararc 360 degrees from 3.3 0.5  center at 3 0.5 
\endpicture} at 0 2
\put{\beginpicture
\put{$N(1)$} at .5 0.5
\put{$1$} at 1 0
\put{$1$} at 2 0
\put{$0$} at 3 0.5
\put{$1$} at 3 0
\put{$1$} at 3 -.5
\arr{1.2 0}{1.8 0}
\arr{2.2 0}{2.8 0}
\arr{2.2 0.1}{2.8 0.4}
\arr{2.2 -.1}{2.8 -.4}
\setdots <.7mm>
\circulararc 360 degrees from 2.3 0  center at 2 0 
\endpicture} at 4 2
\put{\beginpicture
\put{$N(2)$} at .5 0.5
\put{$0$} at 1 0
\put{$1$} at 2 0
\put{$0$} at 3 0.5
\put{$1$} at 3 0
\put{$0$} at 3 -.5
\arr{1.2 0}{1.8 0}
\arr{2.2 0}{2.8 0}
\arr{2.2 0.1}{2.8 0.4}
\arr{2.2 -.1}{2.8 -.4}
\setdots <.7mm>
\circulararc 360 degrees from 2.3 0  center at 2 0 
\endpicture} at 7.5 2
\put{\beginpicture
\put{$N(2)$} at .5 0.5
\put{$0$} at 1 0
\put{$1$} at 2 0
\put{$0$} at 3 0.5
\put{$0$} at 3 0
\put{$1$} at 3 -.5
\arr{1.2 0}{1.8 0}
\arr{2.2 0}{2.8 0}
\arr{2.2 0.1}{2.8 0.4}
\arr{2.2 -.1}{2.8 -.4}
\setdots <.7mm>
\circulararc 360 degrees from 2.3 0  center at 2 0 
\endpicture} at 11 2
\put{\beginpicture
\put{$M$} at .5 0.5
\put{$1$} at 1 0
\put{$3$} at 2 0
\put{$1$} at 3 0.5
\put{$2$} at 3 0
\put{$2$} at 3 -.5
\arr{1.2 0}{1.8 0}
\arr{2.2 0}{2.8 0}
\arr{2.2 0.1}{2.8 0.4}
\arr{2.2 -.1}{2.8 -.4}
\circulararc 360 degrees from 1.3 0  center at 1 0 
\endpicture} at 0 0.2
\put{\beginpicture
\put{$N'(1)$} at .5 0.5
\put{$0$} at 1 0
\put{$1$} at 2 0
\put{$0$} at 3 0.5
\put{$1$} at 3 0
\put{$1$} at 3 -.5
\arr{1.2 0}{1.8 0}
\arr{2.2 0}{2.8 0}
\arr{2.2 0.1}{2.8 0.4}
\arr{2.2 -.1}{2.8 -.4}
\setdots <.7mm>
\circulararc 360 degrees from 2.3 0  center at 2 0 
\endpicture} at 4 0.2
\put{\beginpicture
\put{$N'(2)$} at .5 0.5
\put{$0$} at 1 0
\put{$2$} at 2 0
\put{$1$} at 3 0.5
\put{$1$} at 3 0
\put{$1$} at 3 -.5
\arr{1.2 0}{1.8 0}
\arr{2.2 0}{2.8 0}
\arr{2.2 0.1}{2.8 0.4}
\arr{2.2 -.1}{2.8 -.4}
\setdots <.7mm>
\circulararc 360 degrees from 2.3 0  center at 2 0 
\endpicture} at 7.5 0.2

\endpicture
$$
	\bigskip 
{\bf Example 5.} Radiation modules are usually not exceptional
and, in contrast to exceptional modules, not determined by the dimension vector. 
Here we
exhibit a radiation module $M_1$ (the origin $x$ is encircled)
and the  decomposition of the restriction
$M'_1 = N_1(1)\oplus N_1(2)$ to $Q^x$. This module $M_1$ is an indecomposable
regular module and its dimension vector shows that it has self-extensions.
Note that there are 
corresponding radiation modules $M_2$ and $M_3$ with the same 
dimension vector obtained by permuting the arms on the right. 
$$
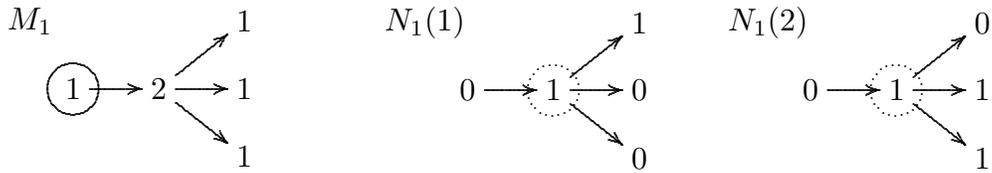
\beginpicture
\setcoordinatesystem units <.95cm,1.5cm>
\put{\beginpicture
\put{$M_1$} at .5 0.5
\put{$1$} at 1 0
\put{$2$} at 2 0
\put{$1$} at 3 0.5
\put{$1$} at 3 0
\put{$1$} at 3 -.5
\arr{1.2 0}{1.8 0}
\arr{2.2 0}{2.8 0}
\arr{2.2 0.1}{2.8 0.4}
\arr{2.2 -.1}{2.8 -.4}
\circulararc 360 degrees from 1.3 0  center at 1 0 
\endpicture} at 0 0
\put{\beginpicture
\put{$N_1(1)$} at .5 0.5
\put{$0$} at 1 0
\put{$1$} at 2 0
\put{$1$} at 3 0.5
\put{$0$} at 3 0
\put{$0$} at 3 -.5
\arr{1.2 0}{1.8 0}
\arr{2.2 0}{2.8 0}
\arr{2.2 0.1}{2.8 0.4}
\arr{2.2 -.1}{2.8 -.4}
\setdots <.7mm>
\circulararc 360 degrees from 2.3 0  center at 2 0 
\endpicture} at 4.5 0
\put{\beginpicture
\put{$N_1(2)$} at .5 0.5
\put{$0$} at 1 0
\put{$1$} at 2 0
\put{$0$} at 3 0.5
\put{$1$} at 3 0
\put{$1$} at 3 -.5
\arr{1.2 0}{1.8 0}
\arr{2.2 0}{2.8 0}
\arr{2.2 0.1}{2.8 0.4}
\arr{2.2 -.1}{2.8 -.4}
\setdots <.7mm>
\circulararc 360 degrees from 2.3 0  center at 2 0 
\endpicture} at 8.5 0
\endpicture
$$
	\bigskip 
We assume now that $(M,x)$ is a radiation module such that $M$ is exceptional.
We want to discuss properties of the corresponding radiation basis of $M$.
Let us start with a general (and well-known) property of exceptional modules:
     \medskip
{\bf Lemma 2.} {\it Let $M$ be a representation of
the quiver $Q$ with $\Ext^1 (M,M) = 0.$ 
Then for every arrow $\alpha$ in $Q$, the
linear map $M_\alpha$ is injective or surjective.}
       \medskip
Proof: Let $\alpha\:x \to y$ be an arrow such that $M_\alpha$ is neither injective nor
surjective.
Write $M_x = M'_x\oplus M''_x$ with $M'_x$ the kernel of $M_\alpha$ and
$M_y = M'_y\oplus M''_y$ with $M'_y$ the image of $M_\alpha$.
By assumption $M'_x$ and $M''_y$ are non-zero, thus there is a non-zero linear map
$N\: M_x \to M_y$ such that $N(M'_x) \subseteq M''_y$ and $N(M''_y) = 0.$

In order to construct a non-trivial self-extension of $M$, start with
the direct sum $M\oplus M$ and replace the map $(M\oplus M)_\alpha =
\left[\smallmatrix M_\alpha & 0 \cr
                      0 & M_\alpha \endsmallmatrix\right]$
by the map
$
\left[\smallmatrix M_\alpha & N \cr
                      0 & M_\alpha \endsmallmatrix\right]$.
The properties of $N$ show that the rank of this map is larger then the rank
of $(M\oplus M)_\alpha$, thus the new representation cannot be isomorphic to
$M\oplus M$, and obviously, it is an extension of $M$ by itself. 
	 \bigskip
{\bf Proposition 3.} {\it Let $(M,x)$ be an exceptional
 radiation module with radiation basis $\Cal B(M,x)$.
Consider an edge
$\{y,z\}$ in $\overline Q$ such that $d(x,y) + 1 = d(y,z).$ 
Then any element of
$\Cal B(M,x)_z$ is connected to precisely one element of $\Cal B(M,x)_y$ by an
arrow.

If $\dim M_y \le \dim M_z$, then any element of $\Cal B(M,x)_y$ is connected to at least one
 element of $\Cal B(M,x)_z$.
If $\dim M_y \ge \dim M_z$, then any element of $\Cal B(M,x)_y$ is connected to at most one
 element of $\Cal B(M,x)_z$.}
	\medskip 
Note that here we consider a path
$$
{\beginpicture
\setcoordinatesystem units <1cm,1cm>
\put{$x$} at 0 0
\put{$y$} at 3 0
\put{$z$} at 4 0
\put{$\alpha$} at 3.5 0.2
\plot 0.2 0  0.8 0 /
\plot 2.2 0  2.8 0 /
\plot 3.2 0  3.8 0 /
\setdots <1mm>
\plot 1.3 0  1.7 0 /
\endpicture} 
$$
in the tree $\overline Q$. 
It follows that the edges between $\Cal B(M,x)_y$ and $\Cal B(M,x)_z$
are of one of the following forms: 
$$
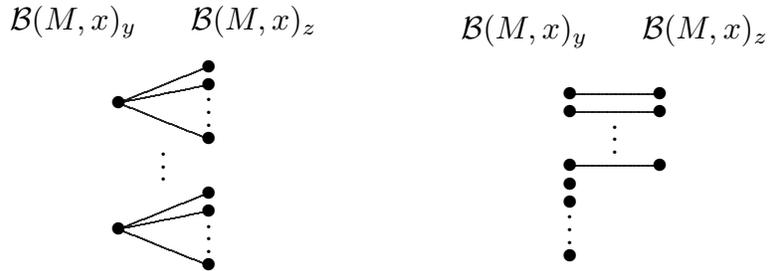
{\beginpicture
\setcoordinatesystem units <1cm,1cm>
\put{\beginpicture
\multiput{$\bullet$} at 0 1.6  1 1.2  1 1.8  1 2 
 0 0.2  1 -.2  1 0.4  1 0.6 /
\multiput{$\vdots$} at   1 0.2  1 1.6  0.5 1 /
\plot 0 0.2  1 -.2 /
\plot 0 0.2  1 0.4 /
\plot 0 0.2  1 0.6 /
\plot 0 1.6  1 1.2 /
\plot 0 1.6  1 1.8 /
\plot 0 1.6  1 2 /

\put{$\Cal B(M,x)_y$} at -.5 2.5
\put{$\Cal B(M,x)_z$} at 1.5 2.5
\endpicture} at 0 0
\put{\beginpicture
\multiput{$\bullet$} at 0 0  0 .6  0 .8  0 1  0 1.6 0 1.8 
   1 1  1 1.6 1 1.8 /
\multiput{$\vdots$} at 0 0.4  0.5 1.4  /
\plot 0 1  1 1 /
\plot 0 1.6  1 1.6 /
\plot 0 1.8  1 1.8 /
\put{} at 0 2
\put{$\Cal B(M,x)_y$} at -.5 2.5
\put{$\Cal B(M,x)_z$} at 1.5 2.5

\endpicture} at 5 0

\endpicture}
$$
	\medskip
{\bf Corollary 1.} {\it Let $(M,x)$ be an exceptional radiation module. Consider an edge
$\{y,z\}$ in $\overline Q$ such that $d(x,y) + 1 = d(y,z).$ Then the number of arrows
between $\Cal B(M,x)_y$ and $\Cal B(M,x)_z$ is precisely $\dim M_z = |\Cal B(M,x)_z|$.}
	\medskip
{\bf Corollary 2.} {\it Let $(M,x)$ be an exceptional radiation module. Let $\alpha$ be an arrow
which connects the vertices $y$ and $z$. If $\dim M_y = \dim M_z,$ then with respect to 
a suitable ordering of 
the bases $\Cal B(M,x)_y$ and $\Cal B(M,x)_z,$ the map $M_\alpha$
is given by the identity matrix.}
     \medskip
Proof of Proposition 3. Write $\Cal B = \Cal B(M,x).$ 
First, we claim that any element of
$\Cal B_z$ is connected to precisely one element of $\Cal B_y$ by an
arrow. This is an immediate consequence of the way the radiation quiver is constructed.
Namely, any vertex $b$ of $R(M,x)$ occurs at a certain step as an origin. In the next step,
this vertex $b$ is connected by an arrow to the new origin, say $b'$.
In later steps, no further arrows are attached which involve $b$.
Thus, let $d(x,z) = i$ and let $D$ be the maximal distance of vertices to $x$.
We start with step $0$, namely with the disjoint union of radiation quivers
consisting of a single vertex, one for each vertex with distance $D$ to $x$.
The basis $\Cal B_z$ is constructed in step $D-i$, all the elements of $\Cal B_z$
are origins of radiation quivers. In the step $D-i+1$, the basis $\Cal B_y$ is
constructed: here, any vertex $b\in \Cal B_z$ is connected by an edge to precisely
one vertex $b'\in \Cal B_y.$ 

For the second assertion, we have to take into account the orientation of the arrow $\alpha$
which connects $y$ and $z$, as well as whether $M_\alpha$  is  injective or surjective. 
	      \medskip 
(1) {\it If $\alpha$ starts at $y$ and $M_\alpha$ is injective, then
any vertex in $\Cal B_y$ is connected by an edge to at least one vertex in $\Cal B_z$.}
Proof: If $b\in \Cal B_y$ is not connected by an edge to a vertex in $\Cal B_z$, 
then $b$ belongs to the kernel of
$M_\alpha$.
	\medskip
(2) {\it If $\alpha$ starts at $y$ and $M_\alpha$ is surjective, then
any vertex in $\Cal B_y$ is connected  by an edge
to at most one vertex in $\Cal B_z$.}
Proof: If there are vertices $b'\neq b''$ in $\Cal B_z$
which are connected by edges to $b\in \Cal B_y$, then
neither $b'$ not $b''$ belong to the image of $M_\alpha$.
	\medskip
(3) {\it If $\alpha$ starts at $z$ and $M_\alpha$ is injective, then
any vertex in $\Cal B_y$ is connected  by an edge
to at most one vertex in $\Cal B_z$.}
Proof: If there are vertices $b'\neq b''$ in $\Cal B_z$
which are connected by edges to $b\in \Cal B_y$, then
$b'-b''$ belongs to the kernel of $M_\alpha$.
	 \medskip
(4) {\it If $\alpha$ starts at $z$ and $M_\alpha$ is surjective, then
any vertex in $\Cal B_y$ is connected  by an edge
to at least one vertex in $\Cal B_z$.}
Proof: If $b\in \Cal B_y$ is not connected by an edge to a vertex in $\Cal B_z$, 
then $b$ is not in the image of $M_\alpha$.
	\bigskip
We should stress that an exceptional radiation module
may have additional tree bases which are not radiation bases.
	\medskip 
{\bf Example 6.}
We deal with the quiver $Q$ of type $\Bbb D_5$ with subspace orientation. 
$$
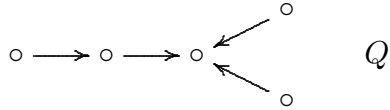
{\beginpicture
\setcoordinatesystem units <1cm,.5cm>
\multiput{$\circ$} at 0 0  1 0  2 0  3 -1  3 1 /
\arr{0.2 0}{0.8 0}
\arr{1.2 0}{1.8 0}
\arr{2.8 -.8}{2.2 -.2}
\arr{2.8 .8}{2.2 .2}
\put{$Q$} at 4 0
\endpicture}
$$
Let $M$ be the maximal indecomposable module. 
There are three vertices 
$x$ with $\dim M_x = 1$, namely the leaves of the quiver $Q$. 
It is easy to see
directly that $(M,x)$ is a radiation module for any leaf $x$ of $Q$,
this is a special case of Proposition 4 shown in the next section. 

Here is a coefficient quiver for $M$, we denote it by $\Gamma$: 
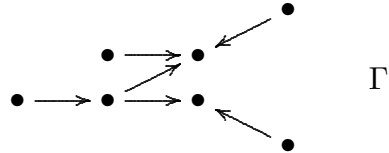
$$
{\beginpicture
\setcoordinatesystem units <1cm,.5cm>
\put{\beginpicture
\multiput{$\bullet$} at 0 0  1 0  1 1  2 0  2 1  3 -1  3 2 /
\arr{0.2 0}{0.8 0}
\arr{1.2 0}{1.8 0}
\arr{1.2 1}{1.8 1}
\arr{2.8 -.8}{2.2 -.2}
\arr{2.8 1.8}{2.2 1.2}
\arr{1.2 0.2}{1.8 0.8}
\put{$\Gamma$} at 4 0.5
\endpicture} at 0 0
\endpicture}
$$
This cannot be a radiation quiver! Let $Q^x$ be obtained by 
removing any leaf $x$ from $Q$, the restriction $M'$ of $M$ to $Q^x$
is the direct sum of two indecomposable modules, thus if $R(M,x)$ is a 
radiation quiver of $M$, then there are two arrows starting in $x$. 
But in $\Gamma$, only one arrow starts at any leaf of $Q$. 
This shows that $\Gamma$ cannot occur as a radiation quiver. 

We also may invoke Corollary 2 above. 
Observe that there is an arrow $\alpha\:y \to z$ such that
$\dim M_y = \dim M_z = 2.$ The matrix presentation given by $\Gamma$ involves three
non-zero coefficients for $M_\alpha$, whereas Corollary 2 assert that the
matrix presentation of $M_\alpha$ with respect to a radiation basis has only two
non-zero coefficients. 
	\bigskip\bigskip
{\bf 4. Dynkin quivers.}
     \bigskip 
{\bf Proposition 4.} {\it Let $Q$ be a Dynkin quiver and $M$ an indecomposable module 
with a thin vertex $x$. 
Then $(M,x)$ is a radiation module.}
     \medskip
This is an immediate consequence of the following lemma.
     \medskip 
{\bf Lemma 3.} 
{\it Let $Q$ be a Dynkin quiver and $M$ an indecomposable module with 
thin vertex $x$. Let $M'$ be the restriction 
of $M$ to the quiver obtained by deleting the vertex $x$. Then $M'$ 
decomposes as the direct sum of orthogonal indecomposable modules $N(i)$.
For every index $i$, there is a unique neighbor $y(i)$ of $x$ with  $\dim N(i)_{y(i)} \neq 0$
and we have $\dim N(i)_{y} = 1$.}
      \medskip
Proof of lemma 3. For a Dynkin quiver, all indecomposable representations are exceptional,
thus we can apply proposition 1(b) and decompose $M' = \bigoplus N(i)^{n(i)}$
where $N(1),\dots, N(t)$ is an  exceptional family. 
For every $i$, there is a unique neighbor $y$ of $x$ with 
$\dim N(i)_{y}\neq 0$. If $\dim N(i)_{y(i)} \ge 2$ 
for some $i$, then the process of simplification (see [R1])
asserts that the full subcategory
of all modules with a filtration with factors $S(x)$ and $N(i)$ is representation
infinite, impossible. This shows that the modules $N(i)$ form an orthogonal  exceptional family.
      \medskip
Note that for $\Delta$ a Dynkin quiver, the only exceptional
modules without a thin vertex are the maximal 
indecomposable modules for the quivers of type $\Bbb E_8$. 
	\bigskip
{\bf Example 7.} Consider the quiver of type $\Bbb E_8$ with subspace orientation and let
$M$ be the maximal indecomposable representation with a thin vertex, say the vertex $x$.
Then $x$ has a unique neighbor $y$ and $\dim M_y = 3.$ It follows that
the restriction $M'$ of $M$ to $Q^x$ is the direct sum of $3$ indecomposable modules
$N(1),N(2),N(3)$. The easiest way to determine these modules $N(i)$ explicitly is as follows:
Consider the hammock $H_y$ for the path algebra $kQ^x,$ this is the support of the
functor $\Hom(P(y),-)$  and has the following shape
$$
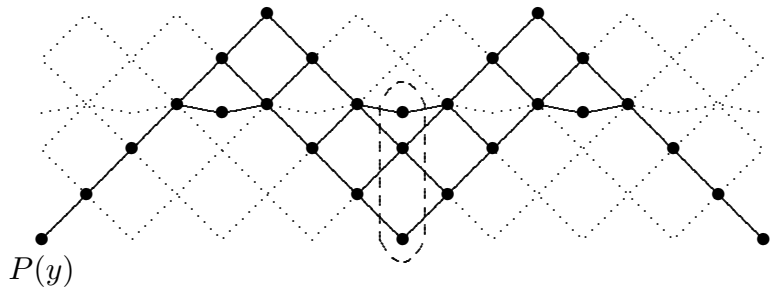
{\beginpicture
\setcoordinatesystem units <.5cm,.5cm>
\setdots <1mm>
\multiput{$\bullet$} at 0 0  1 1  2 2  3 3  4 4  5 5  6 4  7 3  8 2  9 1  10 2  11 3  12 4  13 3
  14 2  15 1  16 0  4 2.8  5 3  6 2  7 1  8 0  8 2.8  9  3  10 4  11 5  12 2.8 /
\plot 0 0  5 5  10 0  15 5  16 4  12 0  7 5  2 0  0 2  3 5  8 0  13 5  16  2  14 0 
    9 5  4 0  0 4  1 5  6 0  11 5  16 0 /
\plot 0 2.8  1 3  2 2.8  3 3  4 2.8  5 3  6 2.8  7 3  8 2.8  9 3  10 2.8  11 3  12 2.8  13 3 
      14 2.8  15 3  16 2.8  / 
\setsolid
\plot 0 0  5 5  9 1  12 4 /
\plot 4 4  8 0  9 1 /
\plot 7 1  11 5  16 0 /
\plot 3 3  4 2.8  5 3  6 4 /
\plot 6 2  7 3  8 2.8  9 3  10 2 /
\plot 10 4  11 3  12 2.8  13 3 /
\setdashes <1mm>
\setquadratic 
\plot 7.5 0  8 -.5  8.5 0  8.5 1.5 
       8.5 3  8 3.5  7.5 3  7.5 1.5  7.5 0 /
\put{} at 0 -.5
\put{$P(y)$} at 0 -.7 
\endpicture}
$$
(it is a poset with partial ordering going from left to right).
There is only one triple of pairwise incomparable elements, it is surrounded
by a dashed line. The corresponding modules $N(1), N(2), N(3)$ 
are the direct summands of $M'$ which we are looking for:
$$
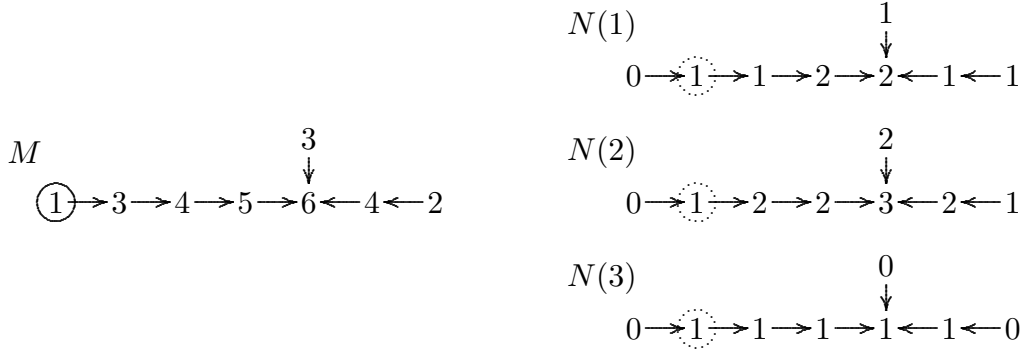
\beginpicture
\setcoordinatesystem units <.7cm,.7cm>
\put{\beginpicture
\put{$M$} at -.5 0.8
\put{$1$} at 0 0
\put{$3$} at 1 0
\put{$4$} at 2 0
\put{$5$} at 3 0
\put{$6$} at 4 0
\put{$4$} at 5 0
\put{$2$} at 6 0
\put{$3$} at 4 1
\arr{0.2 0}{0.8 0}
\arr{1.2 0}{1.8 0}
\arr{2.2 0}{2.8 0}
\arr{3.2 0}{3.8 0}
\arr{4.8 0}{4.2 0}
\arr{5.8 0}{5.2 0}
\arr{4 0.7}{4 0.3}

\circulararc 360 degrees from 0.3 0  center at 0 0 
\endpicture} at -1 0
\put{\beginpicture
\put{$N(1)$} at -.5 0.8
\put{$0$} at 0 0
\put{$1$} at 1 0
\put{$1$} at 2 0
\put{$2$} at 3 0
\put{$2$} at 4 0
\put{$1$} at 5 0
\put{$1$} at 6 0
\put{$1$} at 4 1
\arr{0.2 0}{0.8 0}
\arr{1.2 0}{1.8 0}
\arr{2.2 0}{2.8 0}
\arr{3.2 0}{3.8 0}
\arr{4.8 0}{4.2 0}
\arr{5.8 0}{5.2 0}
\arr{4 0.7}{4 0.3}

\setdots <.7mm>
\circulararc 360 degrees from 1.3 0  center at 1 0 
\endpicture} at 8 2
\put{\beginpicture
\put{$N(2)$} at -.5 0.8
\put{$0$} at 0 0
\put{$1$} at 1 0
\put{$2$} at 2 0
\put{$2$} at 3 0
\put{$3$} at 4 0
\put{$2$} at 5 0
\put{$1$} at 6 0
\put{$2$} at 4 1
\arr{0.2 0}{0.8 0}
\arr{1.2 0}{1.8 0}
\arr{2.2 0}{2.8 0}
\arr{3.2 0}{3.8 0}
\arr{4.8 0}{4.2 0}
\arr{5.8 0}{5.2 0}
\arr{4 0.7}{4 0.3}

\setdots <.7mm>
\circulararc 360 degrees from 1.3 0  center at 1 0 
\endpicture} at 8 0
\put{\beginpicture
\put{$N(3)$} at -.5 0.8
\put{$0$} at 0 0
\put{$1$} at 1 0
\put{$1$} at 2 0
\put{$1$} at 3 0
\put{$1$} at 4 0
\put{$1$} at 5 0
\put{$0$} at 6 0
\put{$0$} at 4 1
\arr{0.2 0}{0.8 0}
\arr{1.2 0}{1.8 0}
\arr{2.2 0}{2.8 0}
\arr{3.2 0}{3.8 0}
\arr{4.8 0}{4.2 0}
\arr{5.8 0}{5.2 0}
\arr{4 0.7}{4 0.3}

\setdots <.7mm>
\circulararc 360 degrees from 1.3 0  center at 1 0 
\endpicture} at 8 -2

\endpicture
$$
	\bigskip
We have mentioned already that an  exceptional representation $M$ of a Dynkin quiver
without a thin vertex is the maximal indecomposable representation of a
quiver $Q$ of type $\Bbb E_8$.
Let us show that the methods presented here provide also for these modules a quite nice
tree basis.

Note that $\bdim M$ looks as follows
$$
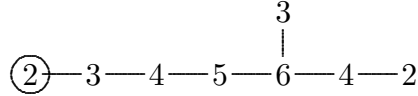
\beginpicture
\setcoordinatesystem units <.7cm,.7cm>
\put{\beginpicture
\put{$2$} at 0 0
\put{$3$} at 1 0
\put{$4$} at 2 0
\put{$5$} at 3 0
\put{$6$} at 4 0
\put{$4$} at 5 0
\put{$2$} at 6 0
\put{$3$} at 4 1
\plot 0.2 0  0.8 0 /
\plot 1.2 0  1.8 0 /
\plot 2.2 0  2.8 0 /
\plot 3.2 0  3.8 0 /
\plot 4.8 0  4.2 0 /
\plot 5.8 0  5.2 0 /
\plot 4 0.7  4 0.3 /

\circulararc 360 degrees from 0.3 0  center at 0 0 
\endpicture} at -1 0

\endpicture
$$
We denote by $x$ the encircled vertex and by $y$ its neighbor.
We consider again the quiver
$Q^x$ obtained from $Q$ by deleting the vertex $x$ (as well as the arrow
involving $x$) and by $M'$ the restriction of $M$ to $Q^x$.
According to Proposition 1(b), we can write 
$M' = \bigoplus_{i=1}^t N(i)^{n(i)}$, where $N(1),\dots, N(t)$ is an exceptional
family of indecomposable modules, all $n(i)\ge 1,$ and all $\dim N(i)_y \ge 1$.

Now $Q^x$ is a quiver of type $\Bbb E_7$, thus any indecomposable representation $N$
of $Q^x$ satisfies $\dim N_y \le 1.$ This shows that $\dim N(i)_y = 1$ for all $i$,
and therefore $t=3$. It follows that $N(1),N(2),N(3)$ is an orthogonal family, 
again we refer to Proposition 1(b). 

As we have seen in Example 7, the hammock $H_y$ for such a quiver of type $\Bbb E_7$
has a unique triple of incomparable elements: the corresponding representations of $Q^x$
are our modules $N(1),N(2),N(3)$. Proposition 4 asserts that the pairs $(N(i),y)$ 
are radiation modules.
We obtain a tree basis of $M$ by using the radiation bases of the pairs $(N(i),y)$
and connecting the origins (they form a basis of $M_y$) with basis elements of $M_x$ as follows:
$$
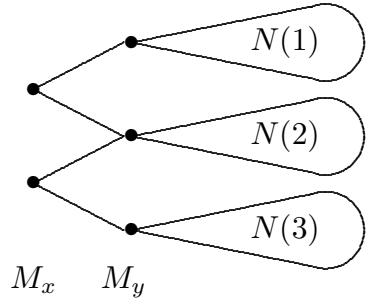
\hbox{\beginpicture
\setcoordinatesystem units <1cm,.8cm>
\put{\beginpicture
\multiput{} at 0 0  2 0.5  2 -.5 /
\plot 0 0  2 0.5 /
\plot 0 0  2 -.5 /
\circulararc -220 degrees from 2 0.5 center at 2.15 0 
\put{$N(2)$} at 1.7 0
\put{$\bullet$} at 0 0 
\endpicture} [l] at 1 0
\put{\beginpicture
\multiput{} at 0 0  2 0.5  2 -.5 /
\plot 0 0  2 0.5 /
\plot 0 0  2 -.5 /
\circulararc -220 degrees from 2 0.5 center at 2.15 0 
\put{$N(3)$} at 1.7 0
\put{$\bullet$} at 0 0 
\endpicture} [l] at 1 -1.3
\put{\beginpicture
\multiput{} at 0 0  2 0.5  2 -.5 /
\plot 0 0  2 0.5 /
\plot 0 0  2 -.5 /
\circulararc -220 degrees from 2 0.5 center at 2.15 0 
\put{$N(1)$} at 1.7 0
\put{$\bullet$} at 0 0 
\endpicture} [l] at 1 1.3
\multiput{$\bullet$} at 0 0.65  0 -.65 / 
\plot 1 1.3  0 0.65  1 0  0 -.65  1 -1.3 /
\put{$M_x\strut$} at 0 -2 
\put{$M_y\strut$} at 1 -2 
\endpicture}
$$
Of course, the dimension vectors of the modules $N(i)$ (in particular the actual radiation
bases of the $N(i)$) depend on the orientation of $Q$. 
	\medskip
{\bf Remark.} In the special cases of dealing with a Dynkin quiver with a unique sink,
Proposition 4 has been obtained independently in the Diplom thesis of V.~Katter,
see [KM].  
	       \bigskip\bigskip
{\bf 5. Bipartite trees without leaves}
     \medskip
Let us now draw the attention to infinite quivers, but we will
assume that the quivers are locally finite without infinite paths.  
For a general discussion of
the Auslander-Reiten components of infinite quiver we may refer to [BLP].

The quiver $Q$ is said to be {\it bipartite} provided that 
every vertex is a sink or a source.  Note that any graph which is a tree can be 
endowed with precisely two orientations so that
we obtain a bipartite tree quiver. 
A vertex $x$ of a tree quiver is called a {\it leaf} provided there 
is at most one arrow $\alpha$ which starts or ends in $x$. 
Note that a tree quiver without leaves has to be infinite. 
As typical examples one should take the $n$-regular tree with bipartite
orientation, where $n\ge 2$. These are the quivers which we later will
use when dealing with the generalized Kronecker quivers. We should stress
that the case $n=3$ was already discussed in [FR]. 

If $Q$ is a quiver, let $Q^*$ be the opposite quiver: it has the same
vertices, but every arrow $\alpha\:x \to y$ is replaced by an
arrow $\alpha^*\:y \to x.$
      \medskip
Let $Q$ be a bipartite tree quiver. Given two vertices $x,y$, let $d(x,y)$ be the distance
between $x$ and $y$.  We denote by $N(x)$ the set of neighbors of $x$ (thus the set
of vertices $y$ with $d(x,y) = 1$).

If $x\in Q_0, t\in \Bbb N_0$, let $B(x,t)$ be the set of
vertices $y$ with $d(x,y) \le t$, we may call it the {\it ball} with center $x$
and radius $t$. If necessary, we will consider $B(x,t)$ as a full subquiver of $Q$.
We will be interested in the balls $B(x,t)$ where $x$ is a sink and
$t$ is even, as well as in those where $x$ is a source and $t$ is odd. 
Note that for these pairs $x,t$, the boundary $\{y\mid d(x,y) = t\}$ of $B(x,y)$
consists of sinks (if $x$ is a sink, then $y$ is also a sink iff
$d(x,y)$ is even. If $x$ is a source, then $y$ is a sink iff $d(x,y)$ is odd).

We denote by $\rho^-$ the composition of all BGP-reflection
functors at all the sources of $Q$, this is a functor $\rho^-\:\mod kQ \to
\mod kQ^*$. Note that $\rho^-\rho^-$ is just the Auslander-Reiten translation
$\tau^-$, see Gabriel [G] (here we use that $\overline Q$ has
no cyclic paths of odd length).

Given a bipartite tree quiver $Q$ without leaves, let us define indecomposable representations 
$P(x,t) = P_Q(x,t)$ 
for certain $x\in Q_0, t\in \Bbb N_0$ as follows: If $x$ is a sink, let
$P(x,0) = S(x) = P(x)$, the one-dimensional representation with support $\{x\}$.
If $x$ is a source, let $P(x,1) = P(x)$, the indecomposable projective module with
top $S(x)$; it is the unique thin module with support $B(x,1)$. 
For $t \ge 2$, define $P(x,t) = \tau^-P(x,t-2)$. Thus, by definition
$$
\align
   P(x,2t) &= \tau^{-t}P(x,0),\cr
   P(x,2t+1) &= \tau^{-t}P(x,1).
\endalign
$$
Looking at the same time at the quivers $Q$ and $Q^*$, we see that for all $t\ge 1$
$$
 P_Q(x,t) = \rho^-P_{Q*}(x,t-1).
$$
We call the modules of the form $P(x,t)$ the {\it preprojective} modules. 
   \medskip
{\bf Lemma 4.} 
\item{(a)}{\it The support of $P(x,t)$ is $B(x,t)$.}

\item{(b)}{\it  If $t\ge 1$ and $d(x,y) \in \{t-1,t\}$, then $\dim P(x,t)_y = 1.$}

\item{(c)}{\it  Let $X$ be an indecomposable representation of $Q$.
Then $\Hom(X,P(x,t)) \neq 0$ if and only if $X$ is preprojective and
the support of $X$ is contained in the support of $P(x,t).$}

\item{(d)} {\it Let $P = P(x,t)$ and $z$ a vertex with $P_z = 0.$ Then
$d(x,z) = t+1$ if and only if $\Ext^1(S(z),P) \neq 0$ and then
$\dim\Ext^1(S(z),P) = \dim P_y = 1.$}
		    \medskip
Proof:
(b) We use induction, the case $t=0$ being trivial. Let $t \ge 1$. We know that
$P_Q(x,t) = \rho^-P_{Q*}(x,t-1)$. Let $y$ be a vertex of $Q$ with $d(x,y) = t.$
There is a unique neighbor $y'$ of $y$ with $d(x,y') = t-1$ and 
$\dim P_{Q*}(x,t-1)_{y'} = 1$ (this is clear for $t=1$, whereas for $t > 1$ this
is the induction hypothesis). The definition of $\rho^-$ shows that 
$P(x,t)_y = P_{Q*}(x,t-1)_{y'}$, thus $\dim P(x,t)_y = 1.$  
If $t\ge 1$ and $y'$ is any vertex of $Q$ with $d(x,y') = t-1$, then $\rho^-$
does not change the vector space at the position $y'$, thus 
$\dim P(x,t)_{y'} = \dim P_{Q*}(x,t-1)_{y'} = 1$.

(a) The induction procedure mentioned in (b) shows that the support of $P(x,t)$
is contained in $B(x,t)$, and as we have see, it contains all the vertices $y$ with 
$d(x,y) = t.$ But the support of an indecomposable module has to be connected,
thus it has to be all of $B(x,t)$.

(c) Let $X$ be an indecomposable representation of $Q$ with $\Hom(X,P(x,t)) \neq 0.$
We write $P(x,t) = \tau^{-s}P$ for some indecomposable projective module (namely
$P = \tau^s P(x,t)$ and $s = \lfloor \frac t2 \rfloor$). If $\tau^s X \neq 0$, then
$\Hom(\tau^s X,P) \neq 0$ and therefore $\tau^s X$ is projective. Altogether, we see
that $X = \tau^{-i}P'$ for some indecomposable projective module $P' = P(y)$ 
and $0 \le i \le s$.
Since $\Hom(P',\tau^{-s+i}P) = \Hom(P',\tau^iP(x,t)) \neq 0$, it follows that
$y$ belongs to the support of $\tau^iP(x,t)$, thus to the support 
$B(x,t)$ of $P(x,t).$ 

Conversely, consider the quiver $B(x,t)$. Since this is a tree quiver, the category of
representations of $B(x,t)$ has a preprojective component. 
Let $\Cal C(x,t)$ be the full subcategory of the indecomposable representations of
$B(x,t)$ which are predecessors of $P(x,t)$, we may consider this as a subcategory
of the category of representations of $Q$. 
As we have seen, any predecessor $X$ of $P(x,t)$ in the category $\mod kQ$
belongs to $\Cal C(x,t).$

(d) Let $P = P(x,t)$ and $z$ a vertex with $P_z = 0$ and $\Ext^1(S(z),P) \neq 0.$
Since $P_z = 0$, we must have $d(x,z) > t.$ Since $\Ext^1(S(z),P) \neq 0$, we 
see that $d(x,z) = t+1$. Then $\dim\Ext^1(S(z),P) = \dim P_y = 1$, according to (b).
    \bigskip
Here is a description of the Auslander-Reiten sequences starting in $P(x,t)$
$$
 0 \to P(x,t) \to \bigoplus\nolimits_{y\in N(x)} P(y,t+1) \to P(x,t+2) \to 0.
$$
	\medskip
We recall the inductive definition of the {\it reachable} objects of a length category:
First of all, the simple projective objects are reachable. Second, if $M$ is
indecomposable, but not simple projective, then $M$ is reachable provided that 
there exists a minimal right almost split map
$M' \to M$ such that all the indecomposable direct summands of $M'$ are reachable.
	\medskip 
{\bf Proposition 5.} {\it The preprojective modules are reachable and they form
a component of the Auslander-Reiten quiver.}
	\medskip
It is of interest to stress the following property of this preprojective component
$\Cal P$: for any indecomposable object $X$ in $\Cal P$, there are sectional paths
starting in a simple projective module $S$ and ending in $X$, namely, 
there is a sectional path from $S = P(y,0)$ to $X = P(x,t)$, provided $d(x,y) = t.$ 
      \medskip
{\bf Proposition 6.} {\it For $t \ge 1$, the pairs $(P(x,t),y)$
with $d(x,y) \in \{t,t-1\}$
are radiation modules.}
    \medskip
Proof. According to Lemma 4(b), we know that the vertices $y$ with $d(x,y) \in \{t,t-1\}$
are thin vertices for $P(x,t)$. For the proof of the proposition, we use induction 
with respect to $t$, the case of $t=1$ being obvious. 

Thus, consider for $t\ge 2$ the pairs $(P(x,t),y)$ with $d(x,y) \in \{t,t-1\}$. The essential
case to deal with is $d(x,y) = t-1$. Namely, if $y$ is a vertex with $d(x,y) = t$,
then there is a unique neighbor $y'$ of $y$ inside the support of $P(x,y)$, and we have both
$\dim P(x,t)_y = 1$ and
$\dim P(x,t)_{y'} = 1.$ Note that $d(x,y') = t-1$. If we know that $(P(x,t),y')$ is a radiation module,
then clearly also $(P(x,t),y)$ is a radiation module.

Let us assume now that $d(x,y) = t-1.$ Let $U$ be the submodule of $P(x,t)$ with
$P(x,t)/U = S(y)$, this is just the restriction of $P(x,t)$ to $Q^{y}$. 
We write $U = \bigoplus_{i\in I} N(i)^{n(i)}$,
with pairwise non-isomorphic indecomposable modules $N(i)$ and integers $n(i) \ge 1$.
Since $U$ is a submodule of $P(x,t)$, we know that these submodules $N(i)$ are again
preprojective modules, see Lemma 4(c). Also, any such module 
$N(i)$ satisfies $N(i)_y = 0$. Of course, we also have 
$\Ext^1(S(y),N(i)) \neq 0,$ and the extensions in $\Ext^1(S(y),N(i))$ are
furnished by an arrow $y \to y(i)$ of $Q$, namely by the unique arrow which connects $y$
with the support of $N(i)$. According to Lemma 4(d), we know that
$\dim N(i)_{y(i)} = 1$ for all $i\in I$. 
Thus, we can apply Proposition 1(b) in order to see 
that the modules
$N(i)$ form an orthogonal exceptional family. It remains to look at the
modules $N(i)$. First of all, if $y(i)$ is a neighbor of $y$ with $d(x,y(i)) = t$,
then $N(i) = S(y(i))$ and the pair $(N(i),y(i))$ is a trivial radiation module.
For the remaining modules $N(i)$ we have $d(x,y(i)) = t-2$. Since
$N(i)$ is preprojective, it is of the form $N(i) = P(x',t')$ for some vertex $x'$ and $t'\le t$.
Since the support of $N(i)$ is properly contained in the support of $P(x,t)$, we even have
$t' < t.$ The vertex
$y(i)$ has to belong to the boundary of the support of $N(i)$. By induction
we know that $(N(i),y(i))$ is a radiation module. Thus we see that $(P(x,t),y)$
is a radiation module. 
	\bigskip
We have seen in the proof of Proposition 6 that for a source $y$ and 
$d(x,y) = t-1$, the kernel
of the canonical epimorphism $P(x,t) \to S(y)$ is a direct sum of orthogonal
bricks. Let us describe this kernel in more details. 
		\medskip
Given a pair $x,y$ of vertices of $Q$, we denote by $[x,y]$ the set of vertices 
lying on the path between $x$ and $y$. If $S$ is a set of vertices of $Q$, we denote by
$N(S)$ the set of neighbors of $S$: these are the vertices $z$ which do not belong to
$S$ but such that there is a vertex $z'\in S$ with $d(z,z') = 1.$
We will be interested in the sets $N([x,y])$, where $x,y$ are vertices with $y$ a source:
$$
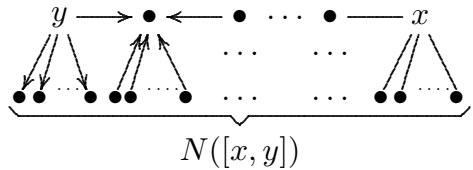
{\beginpicture
\setcoordinatesystem units <1cm,1cm>
\put{} at 0 -1
\put{$y$} at 0 0
\put{$x$} at 4 0
\multiput{$\bullet$} at 1 0  2 0  3 0 /
\arr{0.2 0}{0.8 0}
\arr{1.8 0}{1.2 0}
\plot 3.2 0  3.8 0 /
\multiput{$\hdots$} at 2.5 0  2 -.9  3 -.9  2 -.4  3 -.4 /

\setsolid
\arr{-.1  -.2}{-.4 -.8}
\arr{-.02 -.2}{-.2 -.8}
\arr{.1 -.2}{.3 -.8}

\arr{.65 -.8}{.9  -.2}
\arr{.8 -.8}{.98 -.2}
\arr{1.4 -.8}{1.1 -.2}

\setsolid
\plot 3.9  -.2  3.6 -.8 /
\plot 3.98 -.2  3.8 -.8  /
\plot 4.1 -.2   4.4 -.8 /
\setdots <1mm>
\plot 4. -.8  4.3 -.8 /
\plot 1. -.8  1.3 -.8 /
\plot 0. -.8  0.3 -.8 /
\setsolid

\plot -.55 -1  -.45 -1.1  1.9 -1.1  2 -1.2  2.1 -1.1 4.45 -1.1  4.55 -1 /
\put{$N([x,y])$} at 2 -1.5
\multiput{$\bullet$} at -.44 -.9  -.22 -.9  .35 -.9  
   .62 -.9  .79 -.9  1.4 -.9
      3.56 -.9  3.78 -.9  4.4 -.9 /
      \endpicture}
$$
(This picture reminds on centipedes.)
	  \medskip 
We consider vertices $y\neq z$, where $y$ is a source, and look at the modules
$S(y)$ and $P(z,d(y,z)-1).$
	\medskip

{\bf Lemma 5.} {\it Let $y$ be a source and $z\neq y$ some other vertex. Then
$S(y), P(z,d(y,z)-1)$ is an orthogonal pair with }
$$
\align 
\dim \Ext^1(P(z,d(y,z)-1),S(y)) &= 0\cr
\dim \Ext^1(S(y),P(z,d(y,z)-1)) &= 1.
\endalign
$$
      \medskip
Proof. According to Lemma 4(a), 
we have $P(z,d(y,z)-1)_y = 0$, thus the pair $S(y),P(z,d(y,z)-1)$ is orthogonal. 
We have $\Ext^1(P(z,d(y,z)-1),S(y)) = 0,$ since $S(y)$ is injective.
According to Lemma 4(d), we have $\dim\Ext^1(S(y),P(z,d(y,z)-1)) = 1.$
	  \bigskip
If $d(x,y) = t-1$ and we look at $z\in N([x,y])$ (thus $z$ is one of the legs of the centipedes),
then we have
$$
 d(y,z)-1 = t-d(x,z).
$$
(Namely, let $z'$ be the neighbor of $z$ which belongs to $[x,y]$, then $d(y,z)-1 = d(y,z')$
and $t-d(x,z) = 1+d(x,y)-d(x,z')-1 = d(x,y)-d(x,z') = d(y,z').$) This explains that
we will have to consider modules $P(z,t')$ with $t'= t-d(x,z).$ 
	\medskip 

{\bf Proposition 7.} {\it Let $y$ be a source and $x$ any vertex.
Let $t = d(x,y)+1$. 
Then there is the following exact sequence:
$$
 0 \to \bigoplus_{z\in N([x,y])} P(z,t-d(x,z))
   \to P(x,t) \to S(y) \to 0
$$
The modules $S(y)$ and $P(z,t-d(x,z))$ with $z\in N([x,y])$ are pairwise
orthogonal bricks and the map $P(x,t) \to S(y)$
is the projective cover of $S(y)$ in the full
subcategory }
$$
 \Cal F\Bigl(S(y); P(z,t-d(x,z)), z\in N([x,y])\Bigr).
$$
	\medskip
Proof. We use induction on $t$. In the case $t=1$, we have $x = y$ and there is the
exact sequence
$$ 
  0 \to \bigoplus_{z\in N(x)} P(z,0) \to P(x,1) \to S(x) \to 0.
$$
Now assume the assertion is true for some $t \ge 1.$  Consider a pair of vertices
$x,z_0$ with $d(x,z_0) = t.$ There is a unique vertex $y$ with $d(x,y) = t-1$ and
$d(y,z_0) = 1$. By induction, there is the exact sequence
$$
 0 \to \bigoplus\nolimits_{z\in N([x,y])} P(z,t-d(x,z))
   \to P(x,t) \to S(y) \to 0. \tag{$*$}
$$
Choose some neighbor $z_0$ of
$y$ with $d(x,z_0) = t$. There is (up to isomorphism) a unique module $M$
with top $S(y)$ and socle $S(z_0)$ and we may rearrange the factors in $(*)$ 
in order to obtain an exact sequence of the form
$$
 0 \to \bigoplus\nolimits_{z\in N([x,y]), z\neq z_0} P(z,t-d(x,z))
   \to P(x,t) \to M \to 0.
$$
Now apply the functor $\rho^-$. We obtain the sequence
$$
 0 \to \bigoplus_{z\in N([x,y]),z\neq z_0} P(z,t+1-d(x,z))
   \to P(x,t+1) \to \rho^-M \to 0. \tag{$**$}
$$
The module $\rho^-M$ has top $S(z_0)$ and its socle is the direct sum of the modules
$S(z)$ where $z\neq y$ is a neighbor of $z_0$. Of course, for these vertices $z$,
we have $S(z) = P(z,0)$ and $d(x,z) = t+1$. Since these modules $P(z,0)$ are
projective, we obtain from $(**)$ an exact sequence of the form
$$
 0 \to U \to P(x,t+1) \to S(z_0) \to 0
$$
where $U$ is the direct sum of the modules $P(z,t+1-d(x,z))$ with 
$z\in N([x,y]),z\neq z_0$ as well as those of the form $P(z,0) = P(z,t+1-d(x,z))$ with 
$z\neq y$ a neighbor of $z_0$. It remains to observe that the elements in $N([x,z_0])$
are precisely the vertices $z\in N([x,y]),z\neq z_0$ and the neighbors 
$z\neq y$ of $z_0$. Thus we see that we obtain for the pair $x,z_0$ the required
sequence. 

Now let us show that the modules $S(y)$ and $P(z,t-d(x,z))$ and $z\in N([x,y])$
are orthogonal and that the only non-trivial extension groups between these modules 
are the groups $\Ext^1(S(y),P(z,t-d(x,z))$ and these are one-dimensional.
We can assume that $t\ge 2$ and we denote by $y'$ the unique
vertex with $d(y,y') = 1$ and $d(y',x) = t-2.$
If $d(x,z) = t$, then $P(z,t-d(x,z)) = P(z,0) = S(z)$ is simple projective,
not isomorphic to $S(y)$ and $\dim\Ext^1(S(y),S(z)) = 1.$ Also, for $z\neq z'$,
the modules $S(z), S(z')$ are orthogonal. Thus, assume that $d(x,z) < t.$ 
Let $z'\in [x,y]$ with $d(z,z') = 1.$ Write $a = d(y,z')$ and $b = d(z',x),$
thus $a+b = t-1$ and $d(z,x) = b+1,$ thus $t-d(x,z) = a.$
Since $d(x,z) < t$, we see that $z'\neq y$, thus $d(z,y') = a$ and therefore
$P(z,a)_y = 0$ and $\dim P(z,a)_a = 1$. According to Lemma 4(d), we have
$\dim \Ext^1(S(y),P(z,a)) = 1.$
It follows from $P(z,a)_y = 0$ that 
$P(z,a)_{z'} = 0$ for $d(x,z') = t$, thus $P(z,a)$ and $P(z',0)$ are orthogonal,
and also that $\Ext^1(P(z,a),P(z',0)) = 0.$ Assume now that the we deal with two
vertices $z_1\neq z_2$ in $N([x,y])$ 
with $d(x,z_2)\le d(x,z_1) < t,$ let $a_1 = t-d(x,z_1),$ and
$a_2 = t-d(x,z_2),$  thus $a_1 \le a_2.$ In order to see that there are
no homomorphisms or extensions between $P(z_1,a_1)$ and $P(z_2,a_2)$, 
we can apply $\rho^{a_1}$, thus we have to consider  $P(z_1,0)$
and $P(z_2,a_2-a_1)$. Note that $d(z_1,z_2) = a_2-a_1+2$, thus $z_1$ is
not in the support of $P(z_2,a_2-a_1)$ and not even a neighbor of this support.
This completes the proof.
     \bigskip
{\bf Corollary.} {\it Let $y$ be a source and $x$ any vertex. 
The family of modules $P(z,d(y,z)-1)$ 
with $z \in N([x,y])$ is an exceptional orthogonal family of bricks.}
	\medskip
As an example, let us consider the $3$-regular tree with
bipartite orientation as studied already in [FR]. 
Let us display the dimension vector $\bdim P(x,3)$ for a
source $x$.  As vertex $y$
we take a source with $d(x,y)=2$, thus $(P(x,3),y)$ is a radiation module.
$$
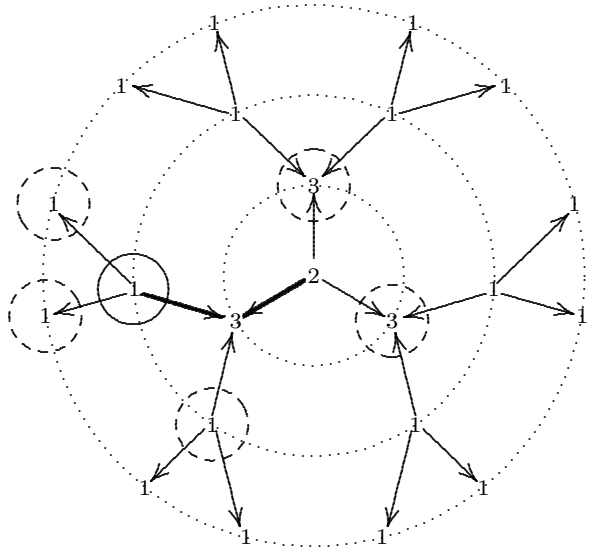
{\beginpicture
\setcoordinatesystem units <1cm,1cm>


\put{$\ssize 2$} at 0 0
\multiput{$\ssize 3$} at -.866 -.5
  .866 -.5  0 1 /
  \setdashes <1mm>
\circulararc 360 degrees from .866 -.9 center at  .866 -.5 
\circulararc 360 degrees from 0 1.4 center at  0 1
\setsolid

\arrow <2mm> [0.25,0.75] from 0 0.2 to  0 0.9
\arrow <2mm> [0.25,0.75] from -0.1 -0.05  to   -.78 -.45
\arrow <2mm> [0.25,0.75] from  0.1 -0.05
 to   .78 -.45
 
\plot -0.1 -0.03   -.78 -.43 /
\plot -0.1 -0.05   -.78 -.45 /
\plot -0.1 -0.07   -.78 -.47 /
\multiput{$\ssize 1$} at -.866 1.8 .866 1.8 /
\arrow <2mm> [0.25,0.75] from  -.78 1.75 to -0.1 1.1
\arrow <2mm> [0.25,0.75] from   .78 1.75 to   0.1 1.1
\multiput{$\ssize 1$} at -2.13 2.1  -1.1 2.8
  2.13 2.1  1.1 2.8 /
  \arrow <2mm> [0.25,0.75] from   -.95 1.8 to  -2. 2.1
  \arrow <2mm> [0.25,0.75] from   -.866 1.9 to -1.07 2.7
  \arrow <2mm> [0.25,0.75] from   .95 1.8 to   2. 2.1
  \arrow <2mm> [0.25,0.75] from  .866 1.9 to  1.07 2.7

  \startrotation by -0.5 0.866 about 0 0
  \multiput{$\ssize 1$} at -.866 1.8 .866 1.8 /
  \arrow <2mm> [0.25,0.75] from -.78 1.75 to  -0.1 1.1
  
\arrow <2mm> [0.25,0.75] from  .78 1.75 to   0.1 1.1

\multiput{$\ssize 1$} at -2.13 2.1  -1.1 2.8
  2.13 2.1  1.1 2.8 /
  \arrow <2mm> [0.25,0.75] from  -.95 1.8  to  -2. 2.1
  \arrow <2mm> [0.25,0.75] from  -.866 1.9 to  -1.07 2.7
  \arrow <2mm> [0.25,0.75] from    .95 1.8 to  2. 2.1
  \arrow <2mm> [0.25,0.75] from    .866 1.9 to 1.07 2.7
  
\setdashes <1mm>
\circulararc 360 degrees from 2.13 2.5 center at  2.13 2.1 
\circulararc 360 degrees from 1.1 2.4 center at  1.1 2.8 

\circulararc 360 degrees from -.866 1.4  center at  -.866 1.8 
\setsolid

\plot .78 1.77  0.1 1.12 /
\plot .78 1.75  0.1 1.1 /
\plot .78 1.73  0.1 1.08 /

\stoprotation

\startrotation by -0.5 -0.866 about 0 0
\multiput{$\ssize 1$} at -.866 1.8 .866 1.8 /
\arrow <2mm> [0.25,0.75] from  -.78 1.75 to -0.1 1.1
\arrow <2mm> [0.25,0.75] from    .78 1.75 to 0.1 1.1
\multiput{$\ssize 1$} at -2.13 2.1  -1.1 2.8
  2.13 2.1  1.1 2.8 /
  \arrow <2mm> [0.25,0.75] from    -.95 1.8 to  -2. 2.1
  \arrow <2mm> [0.25,0.75] from  -.866 1.9 to  -1.07 2.7
  \arrow <2mm> [0.25,0.75] from   .95 1.8 to   2. 2.1
  \arrow <2mm> [0.25,0.75] from    .866 1.9 to 1.07 2.7

  \stoprotation

  \setdots <1mm>
  \circulararc 360 degrees from 1 0 center at 0 0
  \circulararc 360 degrees from 2 0 center at 0 0
  \circulararc 360 degrees from 3 0 center at 0 0
  
\setsolid
\circulararc 360 degrees from -2.3 -.4 center at -2 -.15

\endpicture}
$$
Here, $x$ is the vertex at the center.
The vertex $y$ is encircled using a solid circle, 
the vertices $z\in N([x,y])$ are encircled using dashed circles.
    \bigskip
Here is the radiation quiver $R(P(x,3),y),$ the vertex $y$ has been encircled.
$$
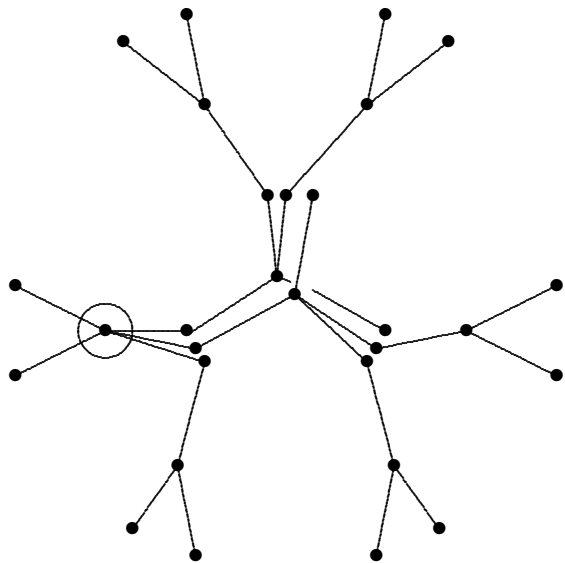
{\beginpicture
\setcoordinatesystem units <1cm,1cm>
\multiput{$\bullet$} at -0.1 0.1  0.1 -0.1 
   -.2 1  0 1  .3 1  -.9 -.85  -1 -.7  -1.1 -.5  .9 -.85  1 -.7  1.1 -.5 
    -1.2 -2  -2 -.5 
   -3 0  -3 -1 
  3 0  3 -1   2 -.5 
   1.2 -2 
   1.7 -2.7  1 -3 
   -1.7 -2.7    -1 -3 
   -.9 2  .9 2 
   -1.1 3  -1.8 2.7   1.1 3  1.8 2.7 /
\plot 0 1  -0.1 0.1  -.2 1 /
\plot .1 -0.1 .3 1 /
\plot -1.1 -.55  -.1 .1 /
\plot -1 -.7  .1 -.1 /
\plot .9 -.85  .1 -.1  1 -.7 /
\plot 1.1 -.5  .3 -.05 /
\plot -.1 .1  0.05 .04 /
\plot  -1.2 -2   -.9 -.85 /
\plot -2 -.5   -.9 -.85 /
\plot -2 -.5  -1 -.7  /
\plot -2 -.5 -1.1 -.5 /
\plot  -3 0  -2 -.5   -3 -1 /
\plot  3 0   2 -.5   3 -1 /
\plot  2 -.5  1 -.7 /
\plot 0.9 -.85  1.2 -2 /
\plot    1.7 -2.7 1.2 -2 1 -3 /
\plot    -1.7 -2.7  -1.2 -2  -1 -3 /
\plot    -1.1 3   -.9 2   -1.8 2.7 /
\plot   1.1 3    .9 2   1.8 2.7 /
\plot   -.9 2  -.2 1 /
\plot   .9 2  0 1 /
\circulararc 360 degrees from -2.3 -.5 center at -2 -.5 

\endpicture}
$$
	\medskip
Finally, let us show the dimension vector 
$\bdim P(x,4)$ with $x$ a sink.
$$
{\beginpicture
\setcoordinatesystem units <1cm,1cm>

\put{$\ssize 7$} at 0 0
\multiput{$\ssize 3$} at -.866 -.5
  .866 -.5  0 1 /
  
\setdashes <1mm>
\circulararc 360 degrees from .866 -.9 center at  .866 -.5 
\circulararc 360 degrees from 0 1.4 center at  0 1
\setsolid

\arrow <2mm> [0.25,0.75] from 0 0.9 to  0 0.2
\arrow <2mm> [0.25,0.75] from -.78 -.45
 to   -0.1 -0.05
 \arrow <2mm> [0.25,0.75] from .78 -.45
  to   0.1 -0.05
  
\plot -0.1 -0.03   -.78 -.43 /
\plot -0.1 -0.05   -.78 -.45 /
\plot -0.1 -0.07   -.78 -.47 /

\multiput{$\ssize 4$} at -.866 1.8 .866 1.8 /
\arrow <2mm> [0.25,0.75] from -0.1 1.1  to  -.78 1.75
\arrow <2mm> [0.25,0.75] from  0.1 1.1  to   .78 1.75
\multiput{$\ssize 1$} at -2.13 2.1  -1.1 2.8
  2.13 2.1  1.1 2.8 /
  \arrow <2mm> [0.25,0.75] from -2. 2.1     to  -.95 1.8
  \arrow <2mm> [0.25,0.75] from -1.07 2.7   to  -.866 1.9
  \arrow <2mm> [0.25,0.75] from 2. 2.1      to  .95 1.8
  \arrow <2mm> [0.25,0.75] from 1.07 2.7    to  .866 1.9

  \multiput{$\ssize 1$} at -3.27 2.3  -2.65 3 -1.7 3.62  -.7 3.93  /
  \multiput{$\ssize 1$} at  3.27 2.3  2.65 3  1.7 3.62  .7 3.93 /

  \arrow <2mm> [0.25,0.75] from -2.3 2.14   to  -3.15 2.3
  \arrow <2mm> [0.25,0.75] from -2.15 2.2    to  -2.6 2.9
  \arrow <2mm> [0.25,0.75] from -1.15 2.9    to  -1.65 3.52
  \arrow <2mm> [0.25,0.75] from -1.05 2.9    to  -.7 3.83

  \arrow <2mm> [0.25,0.75] from 2.3 2.14   to  3.15 2.3
  \arrow <2mm> [0.25,0.75] from 2.15 2.2    to  2.6 2.9
  \arrow <2mm> [0.25,0.75] from 1.15 2.9    to  1.65 3.52
  \arrow <2mm> [0.25,0.75] from 1.05 2.9    to  .7 3.83

  \startrotation by -0.5 0.866 about 0 0
  \multiput{$\ssize 4$} at -.866 1.8 .866 1.8 /
  \arrow <2mm> [0.25,0.75] from -0.1 1.1  to  -.78 1.75
  \arrow <2mm> [0.25,0.75] from  0.1 1.1  to   .78 1.75
  \multiput{$\ssize 1$} at -2.13 2.1  -1.1 2.8
    2.13 2.1  1.1 2.8 /
    \arrow <2mm> [0.25,0.75] from -2. 2.1     to  -.95 1.8
    \arrow <2mm> [0.25,0.75] from -1.07 2.7   to  -.866 1.9
    \arrow <2mm> [0.25,0.75] from 2. 2.1      to  .95 1.8
    \arrow <2mm> [0.25,0.75] from 1.07 2.7    to  .866 1.9

    \multiput{$\ssize 1$} at -3.27 2.3  -2.65 3 -1.7 3.62  -.7 3.93  /
    \multiput{$\ssize 1$} at  3.27 2.3  2.65 3  1.7 3.62  .7 3.93 /

    \arrow <2mm> [0.25,0.75] from -2.3 2.14   to  -3.15 2.3
    \arrow <2mm> [0.25,0.75] from -2.15 2.2    to  -2.6 2.9
    \arrow <2mm> [0.25,0.75] from -1.15 2.9    to  -1.65 3.52
    \arrow <2mm> [0.25,0.75] from -1.05 2.9    to  -.7 3.83

    \arrow <2mm> [0.25,0.75] from 2.3 2.14   to  3.15 2.3
    \arrow <2mm> [0.25,0.75] from 2.15 2.2    to  2.6 2.9
    \arrow <2mm> [0.25,0.75] from 1.15 2.9    to  1.65 3.52
    \arrow <2mm> [0.25,0.75] from 1.05 2.9    to  .7 3.83
    
\setsolid

\plot .78 1.77  0.1 1.12 /
\plot .78 1.75  0.1 1.1 /
\plot .78 1.73  0.1 1.08 /

\plot 2. 2.08   .95 1.78 / 
\plot 2. 2.1   .95 1.8 / 
\plot 2. 2.12   .95 1.82 / 

\circulararc 360 degrees from 2.13 2.5 center at  2.13 2.1 

\setdashes <1mm>

\circulararc 360 degrees from  3.27 2.7 center at   3.27 2.3 
\circulararc 360 degrees from 2.65 3.4  center at  2.65 3 
\circulararc 360 degrees from 1.1 2.4 center at  1.1 2.8 

\circulararc 360 degrees from -.866 1.4  center at  -.866 1.8 
\setsolid

\stoprotation

\startrotation by -0.5 -0.866 about 0 0
\multiput{$\ssize 4$} at -.866 1.8 .866 1.8 /
\arrow <2mm> [0.25,0.75] from -0.1 1.1  to  -.78 1.75
\arrow <2mm> [0.25,0.75] from  0.1 1.1  to   .78 1.75
\multiput{$\ssize 1$} at -2.13 2.1  -1.1 2.8
  2.13 2.1  1.1 2.8 /
  \arrow <2mm> [0.25,0.75] from -2. 2.1     to  -.95 1.8
  \arrow <2mm> [0.25,0.75] from -1.07 2.7   to  -.866 1.9
  \arrow <2mm> [0.25,0.75] from 2. 2.1      to  .95 1.8
  \arrow <2mm> [0.25,0.75] from 1.07 2.7    to  .866 1.9

  \multiput{$\ssize 1$} at -3.27 2.3  -2.65 3 -1.7 3.62  -.7 3.93  /
  \multiput{$\ssize 1$} at  3.27 2.3  2.65 3  1.7 3.62  .7 3.93 /

  \arrow <2mm> [0.25,0.75] from -2.3 2.14   to  -3.15 2.3
  \arrow <2mm> [0.25,0.75] from -2.15 2.2    to  -2.6 2.9
  \arrow <2mm> [0.25,0.75] from -1.15 2.9    to  -1.65 3.52
  \arrow <2mm> [0.25,0.75] from -1.05 2.9    to  -.7 3.83

  \arrow <2mm> [0.25,0.75] from 2.3 2.14   to  3.15 2.3
  \arrow <2mm> [0.25,0.75] from 2.15 2.2    to  2.6 2.9
  \arrow <2mm> [0.25,0.75] from 1.15 2.9    to  1.65 3.52
  \arrow <2mm> [0.25,0.75] from 1.05 2.9    to  .7 3.83
  \stoprotation

  \setdots <1mm>
  \circulararc 360 degrees from 1 0 center at 0 0
  \circulararc 360 degrees from 2 0 center at 0 0
  \circulararc 360 degrees from 3 0 center at 0 0
  \circulararc 360 degrees from 4 0 center at 0 0

  \endpicture}
  $$
  Again, $x$ is the vertex at the center. We have chosen a source
$y$ such that $d(x,y)=3$, the vertex $y$ is encircled
using a solid circle, whereas the dashed circles mark
the vertices $z\in N([x,y])$. Altogether, we see that there is an exact sequence
of the form
$$
 0 \to U \to  
 P(x,4) \to S(y) \to 0
$$
with
$$
 U = P(z_1,0)\oplus P(z_2,0)\oplus P(z_3,1)\oplus P(z_4,2)\oplus P(z_5,3)
 \oplus P(z_6,3)
$$
	\bigskip\bigskip
{\bf 6. The generalized Kronecker quivers and Schofield induction.}
     \medskip
The generalized Kronecker quivers are the quivers $K(n)$ with $2$ vertices, a sink and a source,
and $n$ arrows (going from the source to the sink); the case $n=2$ is the ordinary Kronecker
quiver, its representations have been studied by Weierstra\ss{} and Kronecker. The universal
cover $Q(n)$ of the quiver $K(n)$ is the $n$-regular tree with bipartite orientation. Using 
the push-down functor, any indecomposable representation of the quiver $Q(n)$ yields an indecomposable representation of $K(n).$ In particular, the representation 
$P_{Q(n)}(x,t)$ of $Q(n)$ defined in the last section provides an indecomposable representation
$P_{K(n)}(t)$ and one obtains in this way just the preprojective $K(n)$-modules
(the special case $n=3$ has been discussed in detail in [FR], all the considerations presented there
can easily be adapted to the general case). Obviously, under the push-down functor a tree basis 
is sent to a tree basis. Since all the modules $P(x,t)$ are radiation modules we obtain in this way
distinguished tree bases for the preprojective $K(n)$-modules. The dual considerations yield
distinguished tree bases for the preinjective $K(n)$-modules. 
	      \medskip 
Let us consider now an arbitrary quiver $Q$ and $M$ an exceptional representation of $Q$.
It has been shown in [R3] (see also [R5)] that $M$ is a tree module. Here we want to
outline in which way we can use the previous assertions in order to exhibit a nice
tree basis of $M$. In order to exhibit a tree basis for an exceptional
module $M$, one uses Schofield induction, thus one considers 
exact sequences of the form
$$
 0 \to Y^y \to M \to X^x \to 0
$$
with indecomposable middle term, where $(X,Y)$ is an orthogonal exceptional pair 
with $\dim \Ext^1(X,Y) = e > 0$ and $(x,y)$ is the dimension vector of a sincere
preprojective or preinjective representation $E$ of the $e$-Kronecker quiver $K(e)$.
Note that the triple $(X,Y;E)$ uniquely determines $M$ and one obtains in this way
inductively 
all the exceptional representations of $Q$. In order to construct a tree basis of $M$ we need
to know tree bases of $X,Y$ and $E$. For the procedure to obtain a tree basis of $M$
from the tree bases of $X,Y,E$ we refer to [R3] (see also [R5]).

A required tree basis for $E$ has been exhibited already in [R3], but
the present note provides an intrinsic way for obtaining such a tree basis.
Let us repeat: 
Instead of working with the $e$-Kronecker quiver itself, we consider
its universal cover, the $e$-regular tree $Q(e)$ with bipartite orientation.
Let $E$ be obtained from the representation $\widetilde E$ of $Q(e)$
by the push-down functor. Now  $Q(e)$ is a bipartite tree quiver without
leaves and $\widetilde E$ is a preprojective or preinjective 
representation of $Q$. Thus, there is a vertex $x$
such that the pair $(\widetilde E,x)$ is a radiation quiver.
But this means that $\widetilde E$ has a distinguished tree basis. Under the push-down
functor, we obtain a distinguished tree basis for $E$.
	 \bigskip\bigskip 

{\bf References}
     \medskip
\item{[BLP]} R. Bautista, S. Liu, and C. Paquette, Representation theory of 
  strongly locally finite quivers,   Proc. London Math. Soc.  (to appear).
\item{[FR]} Ph.~Fahr, C.~M.~Ringel:
  Categorification of the Fibonacci numbers using representations of quivers. 
  Journal of Integer Sequences. Vol. 15 (2012), Article 12.2.1
\item{[G]} P.~Gabriel: Auslander-Reiten sequences and representation-finite algebras.
  In Representation theory I. Springer LNM 831, Berlin (1980) 1-71.
\item{[K]} R.~Kinser: Rank functions on rooted tree quivers. Duke Math. J. 152 (2010),
  27-92.
\item{[KM]} V.~Katter, N.~Mahrt: Reduced representations of rooted trees. (To appear).
\item{[R1]} C.~M.~Ringel:   Representations of K--species and bimodules. J. Algebra 41
  (1976), 269-302.
\item{[R2]} C.~M.~Ringel: Tame algebras and integral quadratic forms. Springer LNM 1099 (1984). 
\item{[R3]} C.~M.~Ringel: Exceptional modules are tree modules. 
  Lin. Alg. Appl. 275-276 (1998) 471-493.
\item{[R4]} C.~M.~Ringel: The Gabriel-Roiter measure. 
Bull. Sci. math. 129 (2005), 726-748.
\item{[R5]} C.~M.~Ringel: 
  Indecomposable representations of the Kronecker quivers. 
  Proc. Amer. Math. Soc. 141 (2013), 115-121. 
\item{[RV]} Ringel, C. M., Vossieck, D.: Hammocks. 
   Proc. London Math. Soc. (3)    54 (1987), 216-246.

  \bigskip\bigskip
{\rmk
C. M. Ringel\par
Department of Mathematics, Shanghai Jiao Tong University \par
Shanghai 200240, P. R. China, and \par 
King Abdulaziz University, P O Box 80200\par
Jeddah, Saudi Arabia\par
	\medskip

e-mail: {\ttk ringel\@math.uni-bielefeld.de} \par
}

\bye